\documentclass[preprint,sort&compress,final,10pt]{elsarticle}

\usepackage{amsmath}
\usepackage{amssymb}
\usepackage{graphicx}
\usepackage{latexsym}
\usepackage{epstopdf}
\usepackage{natbib}
\usepackage{color}
\usepackage{hyperref}
\usepackage{amsmath,amscd}

\newtheorem{theorem}{Theorem}[section]

\newtheorem{definition}[theorem]{Definition}
\newtheorem{lemma}[theorem]{Lemma}

\newtheorem{remark}[theorem]{Remark}
\numberwithin{equation}{section}

\def\Proof{\noindent{\bf Proof.}~}
\def\qed{\hfill$\square$\smallskip}

\def\dsum{\displaystyle\sum}

\def\Im{\mathrm{Im}}

\journal{\empty}
\date{}

\begin{document}

\begin{frontmatter}

\title{Persistence of invariant tori in integrable Hamiltonian systems under almost periodic perturbations}

\author[au1]{Peng Huang}

\address[au1]{School of Mathematics Sciences, Beijing Normal University, Beijing 100875, P.R. China.}

\ead[au1]{hp@mail.bnu.edu.cn}

\author[au1]{Xiong Li\footnote{Corresponding author. Partially supported by the NSFC (11571041) and the Fundamental Research Funds for the Central Universities.}}

\ead[au1]{xli@bnu.edu.cn}

\begin{abstract}
In this paper we are concerned with the existence of invariant tori in nearly integrable Hamiltonian systems
\begin{equation*}
H=h(y)+f(x,y,t),
\end{equation*}
where $y\in D\subseteq\mathbb{R}^n$ with $D$ being a closed bounded domain, $x\in \mathbb{T}^n$, $f(x,y,t)$ is a real analytic almost periodic function in $t$ with the frequency ${{\omega}}=(\cdots,{{\omega}}_\lambda,\cdots)_{\lambda\in \mathbb{Z}}\in \mathbb{R}^{\mathbb{Z}}$. As an application, we will prove the existence of almost periodic solutions and the boundedness of all solutions for the second order differential equations with superquadratic potentials depending almost periodically on time.
\end{abstract}


\end{frontmatter}

\section{Introduction}
In this paper we study the existence of invariant tori in the nearly integrable Hamiltonian system
\begin{equation}\label{a8}
H=h(y)+f(x,y,t),
\end{equation}
where $y\in D\subseteq\mathbb{R}^n$ with $D$ being a closed bounded domain, $x\in \mathbb{T}^n$, $f(x,y,t)$ is a real analytic almost periodic function in $t$ with the frequency ${{\omega}}=(\cdots,{{\omega}}_\lambda,\cdots)_{\lambda\in \mathbb{Z}}$ and admits a spatial series expansion
\begin{equation}\label{a7}
f(x,y,t)=\sum \limits _{{k} \in {{\mathbb{Z}}_{\mathcal{S}}^{\mathbb{Z}}}} f_{{k}}(x,y)e^{i \langle {k},{\omega} \rangle t}.
\end{equation}
Here ${\omega}=(\cdots,{\omega}_\lambda,\cdots)_{\lambda\in \mathbb{Z}}$ is a bilateral infinite sequence of rationally independent frequency, that is to say, any finite segments of $\omega$  are rationally independent, \ $\mathcal{S}$ is a family of  finite subsets $A$ of $\mathbb{Z}$, which represents a spatial structure on $\mathbb{Z}$  with $\mathbb{Z}\subseteq \bigcup\limits_{A\in\mathcal{S}} A$  depended on the Fourier exponents $\{\Omega_\lambda: \lambda\in \mathbb{Z}\}$ of a kind of real analytic almost periodic functions  similar to the almost periodic perturbation $f$, and its basis is the  frequency  $\{{\omega}_\lambda: \lambda\in \mathbb{Z}\}$ which is contained in $\{\Omega_\lambda: \lambda\in \mathbb{Z}\}$, ${\mathbb{Z}}_{\mathcal{S}}^{\mathbb{Z}}$ is the space of bilateral infinite integer sequences ${k}=(\cdots,{k}_\lambda,\cdots)_{\lambda\in \mathbb{Z}}$ whose support
$\mbox{supp}\, {k}=\big\{\lambda\ : \ {k}_{\lambda}\neq 0\big\}$ is a finite set of $\mathbb{Z}$ contained in a subset $A$ that belongs to $\mathcal{S}$,  $\langle {k},{\omega} \rangle=\sum \limits _{\lambda\in \mathbb{Z}} {k}_\lambda {\omega}_\lambda$.

Kolmogorov-Arnold-Moser (or KAM) theory is a powerful method about the persistence of quasi-periodic solutions and almost periodic solutions under small perturbations.  KAM theory is not only a collection of specific theorems, but rather a methodology, a collection of ideas of
how to approach certain problems in perturbation theory connected with small divisors.

The classical KAM  theory  was developed for the stability of motions in Hamiltonian systems, that are small perturbations of integrable Hamiltonian systems.  Integrable systems in their phase space contain lots of invariant tori and the classical KAM theory establishes persistence
of such tori,  which carry quasi-periodic motions. The classical KAM  theory concludes that most of invariant tori of integrable Hamiltonian systems can survive uner small perturbation and with Kolmogorov's nondegeneracy condition \cite{{Arnold63},{Eliasson88},{Kolmogorov54},{Poschel01}}.

Later important generalizations of the classical KAM theorem were made
to the R\"{u}ssmann's nondegeneracy condition \cite{{Cheng94},{Russmann90},{Russmann01},{Sevryuk95},{Xu97}}. However, in the case of R\"{u}ssmann's
nondegeneracy condition, one can only get the existence of a family of invariant tori while
there is no information on the persistence of frequency of any torus.

Chow, Li, Yi \cite{Chow02}  and Sevryuk  \cite{Sevryuk06} considered perturbations of moderately degenerate integrable
Hamiltonian systems and proved that the first $d$ frequencies ($d < n$, $n$ denotes the freedom of
Hamiltonian systems) of unperturbed invariant $n$-tori can persist. Xu and You \cite{Xu10} proved that
if some frequency satisfies certain nonresonant condition and topological degree condition,
the perturbed system still has an invariant torus with this frequency under R\"{u}ssmann's
nondegeneracy condition. Zhang and Cheng \cite{Zhang10} concerned with  the persistence of invariant tori for nearly integrable
Hamiltonian systems under time quasi-periodic perturbations, they proved that if the frequency of unperturbed system satisfies the R\"{u}ssmann's nondegeneracy condition and has nonzero Brouwer's topological degree at some Diophantine frequency, then invariant torus with frequency (Diophantine frequency and frequency of time quasi-periodic perturbation) satisfying the Diophantine condition persists under time quasi-periodic perturbations.

However there are only few results available to obtain the existence of almost-periodic solutions
via KAM theory, because it is difficult to treat small divisor problem of infinite frequencies.

In this paper we focus on the almost periodic case, that is, the perturbation in (\ref{a8}) is almost periodic in $t$ and admits the spatial series
expansion, and want to prove most of invariant tori of the integrable Hamiltonian system can survive under small almost periodic perturbations and with Kolmogorov's nondegeneracy condition.

After we get the invariant tori theorem, as an application, we shall prove the existence of almost periodic solutions and the boundedness of all solutions for the differential equation
\begin{equation}\label{a6}
\ddot{x}+x^{2l+1}=\sum_{j=0}^{2l} p_j(t)x^j,\ \ \ x\in\mathbb{R},
\end{equation}
where $p_0,p_1,\cdots,p_{2l}$ are real analytic almost periodic functions with the frequency ${{\omega}}=(\cdots,{{\omega}}_\lambda,\cdots)_{\lambda\in \mathbb{Z}}$.

External forced problem is an important feature of the classical perturbation for Hamiltonian systems.
It is well known that the longtime behaviour of a time dependent nonlinear differential equation
\begin{equation}\label{a1}
\ddot{x}+f(t,x)=0,
\end{equation}
$f$ being periodic in $t$, can be very intricate. For example, there are equations having unbounded solutions but with infinitely many zeros and with nearby unbounded solution having randomly prescribed numbers of zeros and also periodic solution (see \cite{Moser73}).

In contrast to such unboundedness phenomena one may look for conditions on the nonlinearity, in addition to the superlinear condition that
\begin{equation*}
{1\over x}f(t,x)\rightarrow \infty \ \ \ \text{as}\  |x|\rightarrow \infty,
\end{equation*}
which allow to conclude that all solutions of (\ref{a1}) are bounded.

The problem was studied extensively for the differential equation (\ref{a6}) with $p_0,p_1,\cdots,p_{2l}$ being real analytic quasi-periodic functions in $t$ with the frequency ${\omega}=({\omega}_1,{\omega}_2,\cdots,{\omega}_{\mathfrak{m}})$. The first result was  due to Morris \cite{Morris76}, who proved that every solution of equation
\begin{equation}\label{a2}
\ddot{x}+x^3=f(t),
\end{equation}
with $p(t+1)=p(t)$ being continuous, is bounded. This result, prompted by Littlewood in \cite{Littlewood68}, Morris also proved that there are infinitely many quasi-periodic solutions and the boundedness of all solutions of (\ref{a2})  via Moser's twist theorem \cite{Moser62}. This
result was extended to (\ref{a6}) with sufficiently smooth periodic $p_j(t)$ by Dieckerhoff and Zehnder  \cite{Dieckerhoff87}. Later, their result was extended to more general cases by several authors. We refer to \cite{{Levi86},{Laederich91},{Norris92},{You90},{Yuan1}} and references therein.

When $p_0,p_1,\cdots,p_{2l}$ are quasi-periodic, by  using the KAM iterations,\, Levi and Zehnder \cite{Levi95},\, Liu and You \cite{Liu98} independently proved that there are infinitely many quasi-periodic solutions and the boundedness of all solutions for (\ref{a6}) with $p_0,p_1,\cdots,p_{2l}$ being sufficiently smooth and the frequency ${\omega}=({\omega}_1,{\omega}_2,\cdots,{\omega}_{\mathfrak{m}})$ being Diophantine
$$|\langle k,{\omega}\rangle|\geq {{\gamma}\over {|{k}|^\tau}},\ \ \ 0\neq {k}\in \mathbb{Z}^{\mathfrak{m}},$$
for some $\gamma>0,\tau>\mathfrak{m}-1$. On the other hand, by establishing the invariant curve theorem of planar smooth quasi-periodic twist mappings in \cite{Huang0}, we obtained the existence of quasi-periodic solutions and the boundedness of all solutions for an asymmetric oscillation with a quasi-periodic external force (\cite{Huang2016}).

One knows that the Diophantine condition is crucial when applying the KAM theory.  A natural question is whether the boundedness for all solutions, called Lagrangian stability, still holds if ${\omega}=({\omega}_1,{\omega}_2,\cdots,{\omega}_{\mathfrak{m}})$ is not Diophantine but Liouvillean? Wang and You \cite{Wang16} proved the boundedness of all solutions of (\ref{a6}) with $\mathfrak{m}=2$ and ${{\omega}}=(1,\alpha),\alpha\in \mathbb{R}\setminus \mathbb{Q}$, without assuming $\alpha$ to be Diophantine.

Recently,  in \cite{Huang} we established the invariant curve theorem of  planar almost periodic twist mappings, as an application, we proved the existence of almost solutions and the boundedness of all solutions of (\ref{a2}) when $f(t)$ is a real analytic almost periodic function with frequency ${{\omega}}=(\cdots,{{\omega}}_\lambda,\cdots)_{\lambda\in \mathbb{Z}}$. In \cite{Huang1} we also established some variants of the invariant curve theorem on the basis of the invariant curve theorem obtained in \cite{Huang}, and used them to study the existence of almost periodic solutions and the boundedness of all solutions for an asymmetric oscillation with an almost periodic external force.

Before ending the introduction, an outline of this paper is as follows.\ In Section \ref{sec:real},\ we first define real analytic almost periodic functions and their norms, then list some properties of them. The main invariant tori theorem (Theorem \ref{thm3.2}) is given in Section \ref{sec:hal}.\   The proof of the invariant tori theorem  and the measure estimate are given in Sections  \ref{sec:pro}, \ref{sec:KAM}, \ref{sec:ite}, \ref{sec:mea} respectively. In Section \ref{sec:app},\ we will prove the existence of almost periodic solutions and the boundedness of all solutions for (\ref{a6}) with superquadratic potentials depending almost periodically on time.

\section{Real analytic almost periodic functions and their norms}\label{sec:real}
\subsection{The frequency of real analytic almost periodic functions}\label{subsec:frequency}
Throughout the paper we always assume that the real analytic almost periodic function has the Fourier exponents $\{\Omega_\lambda: \lambda\in \mathbb{Z}\}$, and its basis is the frequency $\{{\omega}_\lambda: \lambda\in \mathbb{Z}\}$ which is contained in $\{\Omega_\lambda: \lambda\in \mathbb{Z}\}$. Then for any $\lambda\in \mathbb{Z}$, $\Omega_\lambda$ can be uniquely expressed into
$$
\Omega_\lambda=r_{\lambda_1}{\omega}_{\lambda_1}+\cdots+r_{\lambda_{j(\lambda)}}{\omega}_{\lambda_{j(\lambda)}},
$$
where $r_{\lambda_1},\cdots, r_{\lambda_{j(\lambda)}}$ are rational numbers. Let
$$
\mathcal{S}=\{(\lambda_1,\cdots,\lambda_{j(\lambda)}):\lambda\in\mathbb{Z}\}.
$$
Thus, $\mathcal{S}$ is a family of  finite subsets $A$ of $\mathbb{Z}$, which reflects a spatial structure on $\mathbb{Z}\subseteq \bigcup\limits_{A\in\mathcal{S}} A$.

For the bilateral infinite sequence of rationally independent frequency ${\omega}=(\cdots,{\omega}_\lambda,\cdots)$,
$$\langle {k},{\omega} \rangle=\sum \limits_{\lambda\in \mathbb{Z}} {k}_\lambda {\omega}_\lambda,$$
where due to the spatial structure of $\mathcal{S}$  on $\mathbb{Z}$, ${k}$ runs over integer vectors whose support
$$\mbox{supp}\, {k}=\big\{\lambda\ : \ {k}_{\lambda}\neq 0\big\}$$
is a finite set of $\mathbb{Z}$  contained in a subset $A$ which belongs to $\mathcal{S}$.
Moreover, we define
\begin{equation}\label{a00000}
{\mathbb{Z}}_{\mathcal{S}}^{\mathbb{Z}}:=\Big\{{k}=(\cdots,{k}_\lambda,\cdots)\in\mathbb{Z}^{\mathbb{Z}} :\ \mbox{supp}\,{k}\subseteq A,\ A\in\mathcal{S}\Big\}.
\end{equation}

In the following we will give a norm for real analytic almost periodic functions. Before we describe the norm,  some more definitions  and  notations  are useful.

The main ingredient of our perturbation theory is a nonnegative weight function
$$[\,\cdot\,]\ :\ \ \ A\ \mapsto\ [A]$$
defined on $\mathcal{S}.$ The weight of a subset may reflect its
size,\ its location or something else. Throughout this paper, we always use the following weight function
\begin{equation}\label{a00002}
[A]=1+\sum\limits_{i\in A}log^{\varrho}(1+|i|),
\end{equation}
where $\varrho>2$ is a constant.

In this paper, the frequency ${{\omega}}=(\cdots,{{\omega}}_\lambda,\cdots)$ of  real analytic almost periodic functions is not only rationally independent with $|{\omega}|_{\infty}=\sup\{|{\omega}_\lambda|:\lambda\in \mathbb{Z}\}<+\infty$,  but also satisfies the strongly nonresonant  condition (\ref{b6}) below.

In a crucial way the weight function determines the nonresonance conditions for the small divisors arising in this theory. As we will do later on by means of an appropriate norm,\ it   suffices to estimate these small divisors from below not only in terms of the norm of ${k}$,
$$|{k}|=\sum \limits_{\lambda \in \mathbb{Z}} |{k}_{\lambda}|,$$
but also in terms of the weight of its support
$$[[{k}]]=\min\limits_{\mbox{supp}\, {k} \subseteq A \in \mathcal{S}}[A].$$
Then the nonresonance conditions read
\begin{equation*}
{{| \langle {k},{\omega} \rangle|} \geq {\alpha \over {\Delta\big([[{k}]]\big)\Delta\big(|{k}|\big)}}},\ \ \ \ 0\neq {k} \in \mathbb{Z}_{\mathcal{S}}^{{\mathbb{Z}}},
\end{equation*}
where, as usual, $\alpha$ is a positive parameter and $\Delta$ some fixed approximation function as described in the following. One and the same approximation function is taken here in both places for simplicity, since the generalization is straightforward. A nondecreasing function $\Delta\ :\ [0,\infty) \rightarrow [1,\infty)$ is called an approximation function,\ if
\begin{equation}\label{b15}
 {{\log\Delta(t)}\over t} \searrow 0,\ \ \ \  t \rightarrow\infty,
\end{equation}
and
$$\int_{0}^{\infty} {{\log\Delta(t)}\over {t^2}}\, dt <\infty.$$
In addition,\ the normalization $\Delta(0)=1$ is imposed for definiteness.

In the following we will give a criterion for the existence of strongly nonresonant
frequencies. It is based on growth conditions on the distribution function
$$N_i(t)=\mbox{card}\ \big\{A\in\mathcal{S}\ :\ \mbox{card}(A)=i,\ [A]\leq t\big\}$$
for $i\geq 1$ and $t\geq 0$.

\begin{lemma}\label{lem2.11}
There exist a constant $N_0$ and an approximation function $\Phi$ such that
$$
\begin{array}{ll}
N_i(t)\leq\left\{\begin{array}{ll}
0,\ \ \ \ \ \ \ \ \ \  t<t_i,\\[0.4cm]
N_0 \Phi(t),\ \ \ t\geq t_i
 \end{array}\right.
\end{array}
$$
with a sequence of real numbers $t_i$ satisfying
$$i\log^{\varrho-1} i\leq t_i\sim i\, \log^{\varrho} i$$
for $i$ large with some exponent $\varrho-1>1$. Here we say $a_i \sim b_i$, if there are two constants $c, C$ such that $c a_i\le b_i\le C a_i$ and $c, C$ are independent of $i$.
\end{lemma}
For a rigorous proof of  Lemma \ref{lem2.11} the reader is referred to \cite{Poschel90,Huang}, we omit it here.

According to Lemma \ref{lem2.11}, there exist an approximation function $\Delta$ and a probability measure $\mathfrak{u}$ on the parameter space $\mathbb{R}^{\mathbb{Z}}$ with support at any prescribed point such that the  measure of the set of ${\omega}$ satisfying the following inequalities
\begin{equation}\label{b6}
{{| \langle {k},{\omega} \rangle|} \geq {\alpha \over {\Delta([[{k}]])\Delta(|{k}|)}}},\ \ \ \ \alpha>0,\ \text{for all}\ {k} \in \mathbb{Z}_{\mathcal{S}}^{{\mathbb{Z}}}\backslash\{0\}
\end{equation}
is positive for a suitably small $\alpha$,\ the proof can be found in \cite{Poschel90},\ we omit it here.

Throughout this paper,\ we assume that the frequency  ${{\omega}}=(\cdots,{{\omega}}_\lambda,\cdots)$  satisfying the nonresonance condition (\ref{b6}).

\subsection{The space of  real analytic almost periodic functions}
In order to find almost periodic solutions $x$ for (\ref{a6}), we have to define a kind of real analytic almost periodic functions which admit a spatial  series expansion similar to (\ref{a7}) with  the frequency  ${{\omega}}=(\cdots,{{\omega}}_\lambda,\cdots)$.

We first define  the space of real analytic quasi-periodic functions $Q({\omega})$ as in \cite[chapter 3]{Siegel97}, here the $\mathfrak{m}$-dimensional frequency vector ${\omega}=({\omega}_1,{\omega}_2,\cdots,{\omega}_{\mathfrak{m}})$ is rationally independent, that is, for any $k=(k_1,k_2,\cdots,k_{\mathfrak{m}}) \neq 0$,\ $\langle k,\omega \rangle =\dsum_{j=1}^{\mathfrak{m}} k_j \omega_j \neq 0$.

\begin{definition}\label{def2.1}
A function $f:\mathbb{R} \rightarrow \mathbb{R}$ is called real analytic quasi-periodic with the frequency ${\omega}=({\omega}_1,{\omega}_2,\cdots,{\omega}_{\mathfrak{m}})$,  if there exists a real analytic function
$$F: \theta=(\theta_1,\theta_2,\cdots,\theta_{\mathfrak{m}}) \in \mathbb{R}^{\mathfrak{m}} \rightarrow \mathbb{R}$$
such that $f(t)=F({\omega}_1t, {\omega}_2t,\cdots, {\omega}_{\mathfrak{m}}t)\ \text{for all}\ t\in \mathbb{R}$, where $F$ is $2\pi$-periodic in each variable and bounded in a complex neighborhood $\Pi_{r}^{{\mathfrak{m}}}=\{(\theta_{1},\theta_{2},\cdots,\theta_{{\mathfrak{m}}})\in \mathbb{C}^{\mathfrak{m}} : |\Im\ \theta_{j}|\leq r, j=1,2,\cdots, {\mathfrak{m}} \}$  of\, $\mathbb{R}^{\mathfrak{m}}$ for some $r > 0$. Here we call $F(\theta)$ the shell function of $f(t)$.
\end{definition}

We denote by $Q({\omega})$  the set of real analytic quasi-periodic functions with the frequency ${\omega}=({\omega}_1,{\omega}_2,\cdots,{\omega}_{\mathfrak{m}})$.  Given $f(t)\in Q({\omega})$, the shell function $F(\theta)$ of  $f(t)$ admits a Fourier series expansion
$$F(\theta)=\sum \limits_{{k}\in \mathbb{Z}^{\mathfrak{m}}} f_{{k}}e^{i \langle {k},\theta \rangle },$$
where ${k}=({k}_1,{k}_2,\cdots,{k}_{\mathfrak{m}})$,  ${k}_j$ range over all integers and the coefficients $f_{{k}}$ decay exponentially with $|{k}|=|{k}_1|+|{k}_2|+\cdots+|{k}_{\mathfrak{m}}|$, then $f(t)$ can be represented as a Fourier series of the type from the definition,
$$
\begin{array}{ll}
f(t)=\sum \limits_{{k}\in \mathbb{Z}^{\mathfrak{m}}} f_{{k}}e^{i \langle {k},{\omega} \rangle t}.
\end{array}
$$

In the following we define the norm of the real analytic quasi-periodic function $f(t)$ through that of the corresponding shell function $F(\theta)$.
\begin{definition}
For $r>0$, let $Q_{r}({\omega})\subseteq Q({\omega})$ be the set of real analytic quasi-periodic functions $f$ such that the corresponding shell functions $F$ are bounded on the subset $\Pi_{r}^{{\mathfrak{m}}}$\ with the supremum norm
$$\big|F\big|_{r}=\sup \limits_{\theta\in \Pi_{r}^{{\mathfrak{m}}}}|F(\theta)|=\sup \limits_{\theta\in \Pi_{r}^{{\mathfrak{m}}}}\Big|\sum_{{k}}f_{{k}}e^{i\langle {k},\theta\rangle}\Big|<+\infty.$$
Thus we define $\big|f\big|_{r}:=\big|F\big|_{r}.$
\end{definition}

Similarly, we give the definition of real analytic almost periodic functions with the frequency ${{\omega}}=(\cdots,{{\omega}}_\lambda,\cdots)$ which is not totally arbitrary. Rather, the frequency $\{{\omega}_\lambda: \lambda\in \mathbb{Z}\}$ is a basis contained in the Fourier exponents $\{\Omega_\lambda: \lambda\in \mathbb{Z}\}$, which is given in Subsection \ref{subsec:frequency}. For this purpose, we first define analytic functions on some infinite dimensional space (see \cite{Dineen99}).
\begin{definition}
Let $X$ be a complex Banach space. A function $f : U\subseteq X \rightarrow \mathbb{C}$, where $U$ is an open subset of $X$, is called analytic if $f$ is continuous on $U$, and $f|_{U\cap X_1}$ is analytic in the classical sense as a function of several complex variables for each finite dimensional subspace $X_1$ of $X$.
\end{definition}

\begin{definition}\label{def2.2}
A function $f:\mathbb{R} \rightarrow \mathbb{R}$ is called real analytic almost periodic with the frequency $\omega=(\cdots,\omega_\lambda,\cdots)\in \mathbb{R}^\mathbb{Z}$, if there exists a real analytic function
$$F: \theta=(\cdots,\theta_\lambda,\cdots) \in \mathbb{R}^{\mathbb{Z}} \rightarrow \mathbb{R},$$
which admits a spatial series expansion
$$F(\theta)=\sum \limits _{A\in\mathcal{S}} F_{A}(\theta),$$
where
$$F_{A}(\theta)={\sum \limits _{\mbox{supp}\, k \subseteq A}} f_{k}\,e^{i \langle k,\theta \rangle },$$
such that $f(t)=F(\omega t)\ \text{for all}\ t\in \mathbb{R}$, where $F$ is $2\pi$-periodic in each variable and bounded in a complex neighborhood  $\Pi_{r}=\Big\{\theta=(\cdots,\theta_\lambda,\cdots) \in \mathbb{C}^{\mathbb{Z}} :\ |\Im\,\theta|_{\infty}\leq r\Big\}$ for some $r > 0$, where $|\Im\,\theta|_{\infty}=\sup \limits_{\lambda\in \mathbb{Z}} |\Im\,\theta_{\lambda}|$. Here $F(\theta)$ is called the shell function of $f(t)$.
\end{definition}

Denote by $AP(\omega)$  the set of real analytic almost periodic functions with the frequency $\omega$ defined by Definition \ref{def2.2}. As a consequence of the definitions of $\mathcal{S}$, the support $\mbox{supp}\, k$ of $k$, and $\mathbb{Z}_{\mathcal{S}}^{\mathbb{Z}}$,\ from Definition \ref{def2.2}, the spatial series expansion of the shell function $F(\theta)$ has another form
$$F(\theta)={\sum \limits_{A\in \mathcal{S}}}\ {\sum \limits _{\mbox{supp}\, k \subseteq A}} f_{k}e^{i \langle k,\theta \rangle }=\sum \limits _{k \in {{\mathbb{Z}}_{\mathcal{S}}^{\mathbb{Z}} }} f_{k}e^{i \langle k,\theta \rangle }.$$
Hence $f(t)$ can be represented as a series expansion of the type
\begin{equation}\label{b10004}
f(t)=\sum \limits _{k \in {{\mathbb{Z}}_{\mathcal{S}}^{\mathbb{Z}} }}f_{k}e^{i \langle k,\omega \rangle t}.
\end{equation}
If we define
$$
f_{A}(t)={\sum \limits _{\mbox{supp}\, k \subseteq A}} f_{k}\,e^{i \langle k,\omega \rangle t},
$$
then
$$f(t)=\sum \limits _{A\in\mathcal{S}} f_{A}(t).$$
From the definitions of the support $\mbox{supp}\, k$ and $A$, we know that $f_{A}(t)$ is a real analytic quasi-periodic function with the frequency $\omega_A=\big\{\omega_\lambda\ :\ \lambda\in A\big\}$. Therefore the almost periodic function $f(t)$ can be represented the sum of countably many quasi-periodic functions $f_{A}(t)$ formally.

\subsection{The norms of real analytic almost periodic functions}
Now we can define the norm of the real analytic almost periodic function $f(t)$ through that of the corresponding shell function $F(\theta)$ just like in the quasi-periodic case.

\begin{definition}
Let $AP_{r}(\omega)\subseteq AP(\omega)$ be the set of real analytic almost periodic functions $f$ such that the corresponding shell functions $F$ are real analytic and bounded on the set $\Pi_{r}$\ with the norm
$$\|F\|_{m,r}=\sum \limits _{A\in\mathcal{S}} |F_{A}|_{r}\, e^{m[A]}=\sum \limits _{A\in\mathcal{S}} |f_{A}|_{r}\, e^{m[A]}<+\infty,$$
where $m>0$ is a constant and
$$|F_{A}|_{r}=\sup\limits_{\theta\in\Pi_{r}}|F_{A}(\theta)|=\sup\limits_{\theta\in\Pi_{r}} \Bigg|\sum \limits _{\mbox{supp}\, k \subseteq A}f_{k}\, e^{i\langle k,\theta\rangle}\Bigg|=|f_{A}|_{r}.$$
Hence we define
$$\|f\|_{m,r}:=\|F\|_{m,r}.$$
\end{definition}

\subsection{Properties of  real analytic almost periodic functions}
In the following some properties of real analytic almost periodic functions are given.
\begin{lemma}\label{lem2.4}
The following statements are true:\\
$(i)$ Let $f(t),g(t)\in AP({\omega})$,\ then $f(t)\pm g(t), g(t+f(t))\in AP({\omega});$\\[0.1cm]
$(ii)$ Let $f(t)\in AP({\omega})$ and $\tau=\beta t +f(t)\ (\beta+f'>0,\beta\neq 0)$,\ then the inverse relation is given by $t={\beta^{-1}}\tau +g(\tau)$  and $g\in AP({{\omega} / \beta})$.\ In particular,\ if $\beta=1$,\ then $g\in AP({\omega}).$
\end{lemma}
The detail proofs of Lemma \ref{lem2.4} can be seen in \cite{Huang}, we omit it here.

\section{The Hamiltonian setting and the main result}\label{sec:hal}
Consider the following Hamiltonian
\begin{equation}\label{c1}
H=h(y)+f(x,y,t),
\end{equation}
where $y\in D\subseteq\mathbb{R}^n, x\in \mathbb{T}^n$, $f(x,y,t)$ is a real analytic almost periodic function in $t$ with the frequency ${{\omega}}=(\cdots,{{\omega}}_\lambda,\cdots)$, $D$ is a closed bounded domain.

After introducing two conjugate variables $\theta\in\mathbb{T}^{\mathbb{Z}}$ and $J\in \mathbb{R}^{\mathbb{Z}}$, the Hamiltonian (\ref{c1}) can
be written in the form of an autonomous Hamiltonian as follows
\begin{equation}\label{c2}
H=\langle {\omega},J\rangle+h(y)+F(x,y,\theta),
\end{equation}
where $F(x,y,\theta)$ is the shell function of the almost periodic function $f(x,y,t)$.  Thus, the perturbed motion of Hamiltonian (\ref{c1}) is described by the following equations
\begin{equation}\label{c3}
\begin{aligned}
\dot{\theta} &=H_J={\omega},\\[0.2cm]
\dot{x} &=H_y=h_y(y)+F_y(x,y,\theta),\\[0.2cm]
\dot{J} &=-H_\theta=-F_\theta(x,y,\theta),\\[0.2cm]
\dot{y} &=-H_x=-F_x(x,y,\theta).
\end{aligned}
\end{equation}

When $f=0$, the unperturbed system (\ref{c3})  has invariant tori $\mathcal{T}_0=\mathbb{T}^\mathbb{Z}\times\mathbb{T}^n\times\{0\}\times\{y_0\}$
with the frequency $\overline{{\omega}}=({\omega},\widetilde{{\omega}})$,  carrying an almost periodic flow $\theta={\omega} t, x(t)=x_0+\widetilde{{\omega}}t$, where $\widetilde{{\omega}}=h_y(y_0)$. The aim is to prove the persistence of invariant tori under small perturbations.

We now make the assumption that this system is nondegenerate in the sense that
\begin{equation*}
\text{det}\ h_{yy}=\text{det}\ {{\partial h_y}\over{\partial y}}\neq 0
\end{equation*}
on $D$. Then $h_y$ is an open map, even a local diffeomorphism between $D$ and some open frequency domain $\mathcal{O}\subseteq\mathbb{R}^n$.

As in \cite{Poschel01}, instead of proving  the existence of invariant tori for the Hamiltonian system (\ref{c2}) directly, we are going to concerned with the existence of invariant tori of a family of linear Hamiltonians. This is accomplished by introducing  the frequency
as independent parameters and changing the Hamiltonian system (\ref{c2}) to a parameterized system. This approach was first taken in  \cite{Moser67}.

To this end we write $y=y_0+z$ and expand $h$ around $y_0$ so that
\begin{equation*}
h(y)=h(y_0)+\langle h_y(y_0),z\rangle+\int_{0}^{1}(1-t)\langle h_{yy}(y_t)z,z\rangle dt,
\end{equation*}
where $y_t=y_0+tz$. By assumption, the frequency map is a diffeomorphism
$$h_y : D\rightarrow \mathcal{O},\ \ \ y_0\mapsto \widetilde{\omega}=h_y(y_0).$$
Hence, instead of $y_0\in D$ we may introduce the frequency $\widetilde{\omega}\in \mathcal{O}$ as independent
parameters, determining $y_0$ uniquely. Incidentally, the inverse map is given as
$$g_{\widetilde{{\omega}}} : \mathcal{O}\rightarrow D,\ \ \ \widetilde{{\omega}}\mapsto y_0=g_{{\widetilde{\omega}}}({\widetilde{\omega}}),$$
where $g$ is the Legendre transform of $h$, defined by $g({\widetilde{\omega}})=\sup\limits _{y}(\langle y,{\widetilde{\omega}}\rangle-h(y))$. See
\cite{Arnold78} for more details on Legendre transforms.

Thus we write
\begin{equation*}
H=e({\widetilde{\omega}})+\langle {\omega},J\rangle+\langle {\widetilde{\omega}},z\rangle+F(x,g_{{\widetilde{\omega}}}({\widetilde{\omega}})+z,\theta)+\int_{0}^{1}(1-t)\langle h_{yy}(y_t)z,z\rangle dt
\end{equation*}
and the term $O\big(|z|^2\big)$  can be taken as a new perturbation,
we  obtain the family of Hamiltonians $H=N+P$ with
\begin{equation*}
N=e({\widetilde{\omega}})+\langle {\omega},J\rangle+\langle {\widetilde{\omega}},z\rangle,
\end{equation*}
\begin{equation*}
P(\theta,x,z;{\widetilde{\omega}})=F(x,g_{{\widetilde{\omega}}}({\widetilde{\omega}})+z,\theta)+\int_{0}^{1}(1-t)\langle h_{yy}(g_{{\widetilde{\omega}}}({\widetilde{\omega}})+tz)z,z\rangle dt.
\end{equation*}
They are real analytic in the coordinates $(\theta,x,J,z)$ in $\mathbb{T}^\mathbb{Z}\times\mathbb{T}^n\times\mathbb{R}^\mathbb{Z}\times B$, $B$ some sufficiently
small ball around the origin in $\mathbb{R}^n$, as well as the frequency ${\widetilde{\omega}}$ taken from a parameters domain $\mathcal{O}$ in $\mathbb{R}^n$.

From the definition of shell function of almost periodic  function, $P$ admits a spatial series expansion
$$P(\theta,x,z;{\widetilde{\omega}})=\sum \limits_{A\in\mathcal{S}}P_A(\theta,x,z;{\widetilde{\omega}})=\sum \limits_{A\in\mathcal{S}}\sum \limits _{\  \mbox{supp}\, {k} \subseteq A }P_{A,{k}}(x,z;\widetilde{\omega})\,e^{i\langle {k},\theta\rangle}.$$
Expanding $P_{A,{k}}(x,z;\widetilde{\omega})$ into a Fourier series at $x\in \mathbb{T}^n$, therefore $P$ admits a spatial series
\begin{equation*}
\begin{aligned}
P(\theta,x,z;{\widetilde{\omega}})&=\sum \limits_{A\in\mathcal{S}}P_A(\theta,x,z;{\widetilde{\omega}})\\[0.2cm]
&=\sum \limits_{A\in\mathcal{S}}\sum \limits _{\substack{\mbox{supp}\, {k} \subseteq A   \\  \widetilde{k}\in \mathbb{Z}^n }}P_{A,{k},{\widetilde{k}}}(z;\widetilde{\omega})e^{i(\langle {k},\theta\rangle+\langle {\widetilde{k}},x\rangle)},
\end{aligned}
\end{equation*}
where $\widetilde{k}=(\widetilde{k}_1,\widetilde{k}_2,\cdots,\widetilde{k}_n)$. Let $B=\{1,2,\cdots,n\}$ be a  finite subset of  $ \mathbb{Z}$, and  $\mbox{supp}\, \widetilde{k}\subseteq B$ for all $\widetilde{k}\in \mathbb{Z}^n$.  Denote $\widetilde{\mathcal{S}}=\{A\times B : A\in\mathcal{S}\}$, then $P$ can be represented as a  spatial series of the type
\begin{equation}\label{d103}
P=\sum \limits_{\widetilde{A}\in\widetilde{\mathcal{S}}}P_{\widetilde{A}}=\sum \limits_{\widetilde{A}\in\widetilde{\mathcal{S}}}\sum \limits _{\mbox{supp}\, (k,\widetilde{k}) \subseteq \widetilde{A} }P_{\widetilde{A},{k},{\widetilde{k}}}(z;\widetilde{\omega})e^{i(\langle k,x\rangle+\langle \widetilde{k},\theta\rangle)},
\end{equation}
where $\mbox{supp}\, (k,\widetilde{k})=\mbox{supp}\, {k}\times \mbox{supp}\, \widetilde{k}, P_{\widetilde{A},{k},{\widetilde{k}}}=P_{A,{k},{\widetilde{k}}}(z;\widetilde{\omega})$. Moreover, we define
\begin{equation*}
{\mathbb{Z}}_{\mathcal{\widetilde{S}}}^{\mathbb{Z}}:=\Big\{(k,\widetilde{k})\in\mathbb{Z}^{\mathbb{Z}}\times \mathbb{Z}^n:\ \mbox{supp}\, (k,\widetilde{k}) \subseteq \widetilde{A} ,\ \widetilde{A}\in\widetilde{\mathcal{S}}\Big\}
\end{equation*}
and
$$[[(k,\widetilde{k})]]=\min\limits_{\mbox{supp}\,  (k,\widetilde{k}) \subseteq \widetilde{A} \in \widetilde{\mathcal{S}}}[\widetilde{A}].$$

To state the invariant tori theorem, we therefore singer out the subsets $\mathcal{O}_{\alpha}\subseteq \mathcal{O}$ denote the set of all ${\widetilde{\omega}}$ satisfying
\begin{equation}\label{c000000004}
{{|\langle {k},{\omega} \rangle+ \langle {\widetilde{k}},{\widetilde{\omega}}\rangle|} \geq {{\alpha} \over {\Delta\big([[(k,\widetilde{k})]]\big)\Delta\big(|k|+|\widetilde{k}|\big)}}},\ \ \ \ \text{for all}\ (k,{\widetilde{k}}) \in \mathbb{Z}_{\widetilde{\mathcal{S}}}^{{\mathbb{Z}}}\backslash \{\widetilde{k}=0\},
\end{equation}
with fixed ${{\omega}}$, where $\Delta$ is an approximation function.

\begin{remark}
For the frequency ${\omega}$ being fixed, there is an approximation function $\Delta$ such that the set $\mathcal{O}_{\alpha}$ is a set
of positive Lebesgue measure provided that $\alpha$ is small (See Theorem \ref{thm2.12} in Section \ref{sec:mea}).
\end{remark}

\begin{remark}
From the proofs of Theorem \ref{thm2.12} and the  measure of the set of ${\omega}$ in \cite{Poschel90}, we know that there is the same approximation function $\Delta$ such that (\ref{b6})   and (\ref{c000000004}) hold simultaneously.
\end{remark}

To state the basic result quantitatively we need to introduce a few notations. Let
\begin{equation*}
\begin{aligned}
\mathcal{D}_{r,s}= &\big\{(\theta,x,J,z)\, :\,   |\Im\ \theta|_{\infty}< r, |\Im\ x|< r, |J|_w< s,  |z|< s\big\}\\[0.2cm]
&\subseteq \mathbb{C}^\mathbb{Z}/2\pi\mathbb{Z}^\mathbb{Z}\times\mathbb{C}^n/2\pi\mathbb{Z}^n\times\mathbb{C}^\mathbb{Z}\times\mathbb{C}^n
\end{aligned}
\end{equation*}
and
$$\mathcal{O}_h=\{{\widetilde{\omega}} : |{\widetilde{\omega}}-\mathcal{O}_{\alpha}|<h\}\subset \mathbb{C}^n,$$
where
$$|\theta|_{\infty}=\sup \limits_{\lambda\in \mathbb{Z}} |\theta_{\lambda}|,\ \ \ |J|_w=\sum \limits_{\lambda\in \mathbb{Z}} |J_\lambda|e^{w[\lambda]},$$
$|\cdot|$ stands for the sup-norm of real vectors respectively, where $w\geq 0$  is another parameter, and the weights at the individual lattice sites are defined by
$[\lambda]=\min \limits _{\lambda\in A\in \mathcal{S}}[A].$

Its size is measured in terms of
the weighted norm
$$|||P|||_{m,r,s,h}=\sum \limits _{\widetilde{A}\in\widetilde{\mathcal{S}}} \|P_{\widetilde{A}}\|_{r,s,h}\, e^{m[\widetilde{A}]},$$
where
$$\|P_{\widetilde{A}}\|_{r,s,h}=\sum \limits _{\mbox{supp}\, (k,\widetilde{k}) \subseteq \widetilde{A} }\big|P_{\widetilde{A},k,\widetilde{k}}\big|_{s,h}\,e^{r(|k|+|\widetilde{k}|)},$$
the norm $|\cdot|_{s,h}$ is the sup-norm over $|z|< s$ and $\widetilde{\omega}\in\mathcal{O}_h$.

Now, the Hamiltonians
\begin{equation}\label{c4}
H(\theta,x,J,z;{\widetilde{\omega}})=N+P
=e({\widetilde{\omega}})+\langle {\omega},J\rangle +\langle {\widetilde{\omega}},z\rangle+P(\theta,x,z;{\widetilde{\omega}}),
\end{equation}
$P$ is real analytic on $\mathcal{D}_{r,s}\times\mathcal{O}_h$. The corresponding
Hamiltonian system (\ref{c3}) becomes
\begin{equation*}
\begin{aligned}
\dot{\theta} &=H_J={\omega},\\[0.2cm]
\dot{x} &=H_z={\widetilde{\omega}}+P_z(\theta,x,z;{\widetilde{\omega}}),\\[0.2cm]
\dot{J} &=-H_\theta=-P_\theta(\theta,x,z;{\widetilde{\omega}}),\\[0.2cm]
\dot{z} &=-H_x=-P_x(\theta,x,z;{\widetilde{\omega}}).
\end{aligned}
\end{equation*}
Thus, the persistence of invariant tori for nearly integrable Hamiltonian system (\ref{c3}) is
reduced to the persistence of invariant tori for the family of Hamiltonian systems (\ref{c4})
depending on the parameter ${\widetilde{\omega}}$.  Our aim is to prove the persistence of
the invariant torus
$$
\mathcal{T}_0=\mathbb{T}^\mathbb{Z}\times\mathbb{T}^n\times\{0,0\}
$$
of maximal dimension together with its constant vector field $(\omega,{\widetilde{\omega}})$.

The smallness condition of the following theorem is expressed in terms of two
functions $\Psi_0,\Psi_1$ that are defined on the positive real axis entirely in terms of the
approximation function $\Delta$ and reflect the effect of the small divisors in solving the
nonlinear problem. See Appendix A in \cite{Poschel90} for their definition.

Now we are in a position to  state our main result.

\begin{theorem}\label{thm3.2}
Suppose that $P$ admits a spatial expansion as in (\ref{d103}), is real analytic on  $\mathcal{D}_{r,s}\times\mathcal{O}_h$ and satisfies
the estimate
$$s^{-1}|||P|||_{m,r,s,h}\leq{{\alpha \varepsilon_*}\over \Psi_0(\mu)\Psi_1(\rho)}\leq {h\over{2^6}}$$
for some $0 < \mu \leq m-w$ and $0<\rho< {r/2}$, where $\varepsilon_*=2^{-22}$ is an absolute positive
constant, $\Psi_0(\mu),\Psi_1(\rho)$ is defined by (\ref{f10001}). Then there exists a transformation
$$\mathcal{F} : \mathcal{D}_{r-2\rho,{s/ 2}}\times\mathcal{O}_{\alpha}\rightarrow \mathcal{D}_{r,s}\times\mathcal{O}_h$$
that is real analytic and symplectic for each ${\widetilde{\omega}}$ and uniformly continuous in ${\widetilde{\omega}}$, such that
$$(N+P)\circ\mathcal{F}=e_*+\langle {\omega},J\rangle+\langle {\widetilde{\omega}},z\rangle+\cdots,$$
where the dots denote terms of higher order in $z$. Consequently, the perturbed
system has a real analytic invariant torus of maximal dimension and with a vector field
conjugate to $({\omega},{\widetilde{\omega}})$ for each frequency vector ${\widetilde{\omega}}$ in $\mathcal{O}_{\alpha}$.  These tori are close of order
$s^{-1}|||P|||_{m,r,s,h}$ to the torus $\mathcal{T}_0$ with respect to the norm $|\cdot|_w$.
\end{theorem}

\section{Outline of the Proof of Theorem \ref{thm3.2}}\label{sec:pro}
Theorem \ref{thm3.2} is proven by the familiar KAM-method employing a rapidly converging iteration scheme \cite{Arnold63,Kolmogorov54,Moser67}.
At each step of the scheme, a Hamiltonian
$$H_j=N_j+P_j$$
is considered, which is a small perturbation of some normal form $N_j$.   A transformation $\mathcal{F}_j$ is set up so that
$$H_j\circ \mathcal{F}_j=N_{j+1}+P_{j+1}$$
with another normal form $N_{j+1}$ and a much smaller error term $P_{j+1}$. For instance,
$$|||P_{j+1}|||\leq C_j|||P_n|||^{\kappa}$$
for some $\kappa>1$.  This transformation consists of a symplectic change of coordinates
$\Phi_j$ and a subsequent change $\varphi_j$ of the parameters ${\widetilde{\omega}}$ and is found by linearising
the above equation. Repetition of this process leads to a sequence of transformations
$\mathcal{F}_0,\mathcal{F}_1,\cdots$,  whose infinite product transforms the initial Hamiltonian $H_0$ into a normal
form $N_*$ up to first order.

Here is a more detailed description of this construction. To describe one cycle of this iterative scheme in more detail we now drop
the index  $j$.

Approximating the
perturbation $P$ in a suitable way we write
\begin{equation*}
\begin{aligned}
H&=N+P\\
&=N+R+(P-R).
\end{aligned}
\end{equation*}
In particular, $R$ is chosen such that its
spatial series expansion is finite, hence all subsequent operations are finite dimensional.

The coordinate transformation $\Phi$ is written as the time-1-map of the flow  $X_F^t$
of a Hamiltonian vector field $X_F$ :
$$\Phi=X_F^t\big|_{t=1}.$$
This makes $\Phi$ symplectic. Moreover, we may expand $H\circ\Phi=H\circ X_F^t\big|_{t=1}$
with respect to $t$ at 0 using Taylor's formula. Recall that
$${{d}\over{dt}}G\circ X_F^t=\{G,F\}\circ X_F^t,$$
the Poisson bracket of $G$ and $F$ evaluated at $X_F^t$. Thus we may write
\begin{equation*}
\begin{aligned}
(N+R)\circ\Phi&=N\circ X_F^t\big|_{t=1}+R\circ X_F^t\big|_{t=1}\\
&=N+\{N,F\}+\int_0^1(1-t) \{\{N,F\},F\}\circ X_F^t\, dt\\
&\ \ \ \ \ \  \  +R+\int_0^1\{R,F\}\circ X_F^t\, dt\\
&=N+R+\{N,F\}+\int_0^1\{(1-t) \{N,F\}+R,F\}\circ X_F^t\, dt.
\end{aligned}
\end{equation*}
The last integral is of quadratic order in $R$ and $F$ and will be part of the new error
term.

The point is to find $F$ such that $N + R + \{N, F\} = N_+$ is a normal form.
Equivalently, setting $N_+=N+\widehat{N}$, the linear equation
\begin{equation}\label{d100000001}
\{F,N\}+\widehat{N}=R
\end{equation}
has to be solved for $F$ and $\widehat{N}$, when $R$ is given. Given such a solution, we obtain
$(1 - t)\{N, F\} + R = (1- t)\widehat{N} + t R$ and hence  $ H\circ \Phi = N_+ + P_+$ with
$$P_+=\int_0^1\{(1-t)\widehat{N}+tR,F\}\circ X_F^t\, dt+(P-R)\circ\Phi.$$

Setting up the spatial expansions for $F$ and $\widehat{N}$ of the same form as that for $R$, the
linearized equation (\ref{d100000001}) breaks up into the component equations
$$i\big(\langle {k},{\omega}\rangle+\langle  {\widetilde{k}},{\widetilde{\omega}}\rangle\big)F_{\widetilde{A}}+\widehat{N}_{\widetilde{A}}=R_{\widetilde{A}}.$$
Their solution is well-known and straightforward. These equations introduce the small divisor,  which in our case are zero if and
only if $(k,{\widetilde{k}})$ is zero by the nonresonance conditions. It therefore suffices to choose
$$\widehat{N}_{\widetilde{A}}=[R_{\widetilde{A}}]$$
the mean value of $R_{\widetilde{A}}$ over $\mathbb{T}^{\widetilde{A}}$, and to solve uniquely
$$i\big(\langle {k},{\omega}\rangle+\langle  {\widetilde{k}},{\widetilde{\omega}}\rangle\big) F_{\widetilde{A}}=R_{\widetilde{A}}-[R_{\widetilde{A}}],\ \ \  \ [F_{\widetilde{A}}]=0.$$
We obtain
\begin{equation}\label{100001}
F_{\widetilde{A}}=\sum \limits _{\substack{\mbox{supp}\,(k,\widetilde{k})\subseteq \widetilde{A}\\ (k,\widetilde{k})\neq 0}}{{R_{{\widetilde{A}},k,\widetilde{k}} }\over {i\big(\langle {k},{\omega}\rangle+\langle  {\widetilde{k}},{\widetilde{\omega}}\rangle\big)}}e^{i(\langle {k},\theta\rangle+\langle {\widetilde{k}},x\rangle)},
\end{equation}
where $R_{{\widetilde{A}},k,\widetilde{k}}$ are the Fourier coefficients of $R_{\widetilde{A}}$.

The truncation of $P$ will be chosen so that $R$ is independent of $J$ and  of first order in $z$. Hence the
same is true of each of the $\widehat{N}_{\widetilde{A}}$ and so
$$\widehat{N}=\sum \limits _{{\widetilde{A}}\in \mathcal{S}}\widehat{N}_{\widetilde{A}}=\hat{e}+\langle v({\widetilde{\omega}}),z\rangle.$$
It suffices to change parameters by setting
\begin{equation}\label{100000}
{\widetilde{\omega}}_+={\widetilde{\omega}}+ v({\widetilde{\omega}})
\end{equation}
to obtain a new normal form $N_+ = N +\widehat{N}$. This completes one cycle of the iteration.

By the same truncation, $F$ is  independent of $J$ and of first order in $z$. It follows that $\Phi=X_F^t\big|_{t=1}$  has the form
\begin{equation*}
\begin{array}{ll}
\begin{array}{ll}
\theta=\theta_+,\\[0.2cm]
x=U_1(\theta_+,x_+),\\[0.2cm]
J=U_2(\theta_+,x_+)+U_3(\theta_+,x_+)z_+,\\[0.2cm]
z=U_4(\theta_+,x_+)+U_5(\theta_+,x_+)z_+,
\end{array}
\end{array}
\end{equation*}
where the dependence of all coefficients on ${\widetilde{\omega}}$ has been suppressed. This map is
composed with the inverse $\varphi$ of the parameter map (\ref{100000}) to obtain $\mathcal{F}$.

Such symplectic transformations form a group under composition. So, if
$\mathcal{F}_0,\mathcal{F}_1,\cdots,\mathcal{F}_{j-1}$ belong to this group, then so does $\mathcal{F}^j = \mathcal{F}_0\circ\mathcal{F}_1\circ\cdots\circ\mathcal{F}_{j-1}$ and the
limit transformation $\mathcal{F}$ for $j\rightarrow\infty$.

\section{The KAM step}\label{sec:KAM}
Before plunging into the details of the KAM-construction we observe that it
suffices to consider some normalized value of $\alpha$, say
$$\widetilde{\alpha}=2$$
Indeed, stretching the time scale by the factor $2/\alpha$ the Hamiltonians $H$ and $N$ are
scaled by the same amount, and so is the frequency $\widetilde{\omega}$. By a similar scaling of the
action-variables $J,z$ the radius $s$ may also be normalized to some convenient value.
We will not do this here.

\subsection{The set up}
Consider a Hamiltonian of the form
\begin{equation*}
H=N+P=e({\widetilde{\omega}})+\langle {\omega},J\rangle+\langle {\widetilde{\omega}},z\rangle+P(\theta,x,z;{\widetilde{\omega}}),
\end{equation*}
where
\begin{equation*}
P=\sum \limits_{\widetilde{A}\in\widetilde{\mathcal{S}}}P_{\widetilde{A}}=\sum \limits_{\widetilde{A}\in\widetilde{\mathcal{S}}}\sum \limits _{\mbox{supp}\, (k,\widetilde{k}) \subseteq \widetilde{A} }P_{\widetilde{A},{k},{\widetilde{k}}}e^{i(\langle k,x\rangle+\langle \widetilde{k},\theta\rangle)}.
\end{equation*}
Assume that $P$ is real analytic on the complex domain
$$\mathcal{D}_{r,s}\times\mathcal{O}_h : \ \ |\Im\ \theta|_{\infty}< r, |\Im\ x|< r, |J|_w< s, |z|< s,  |{\widetilde{\omega}}-\mathcal{O}_{\ast}|<h,$$
where $\mathcal{O}_{\ast}$ is a closed subset of the parameter space $\mathbb{R}^n$ consisting of the frequency ${\widetilde{\omega}}$ that
satisfying
\begin{equation}\label{d1}
{{|\langle {k},{\omega} \rangle+ \langle {\widetilde{k}},{\widetilde{\omega}}\rangle|} \geq {{\widetilde{\alpha}} \over {\Delta\big([[(k,\widetilde{k})]]\big)\Delta\big(|k|+|\widetilde{k}|\big)}}},\ \ \ \ \text{for all}\ (k,{\widetilde{k}}) \in \mathbb{Z}_{\widetilde{\mathcal{S}}}^{{\mathbb{Z}}}\backslash \{\widetilde{k}=0\},
\end{equation}
where $\widetilde{\alpha}=2$. Moreover, assume that for some $m>w$,
\begin{equation}\label{d01}
|||H-N|||_{m,r,s,h}=|||P|||_{m,r,s,h}\leq \varepsilon
\end{equation}
is sufficiently small. The precise condition will be given later in the course of the
iteration.

Unless stated otherwise the following estimates are uniform with respect to ${\widetilde{\omega}}$.
Therefore the index $h$ is usually dropped.

\subsection{Truncating the perturbation} \label{subsec:truncating}
Let $\mu$ and $\rho$ be two small and $K$ a large positive parameter to be chosen during
the iteration process. The Fourier series of the $\widetilde{A}$-component $P_{\widetilde{A}}$ of the perturbation
is truncated at order $(|k|+|\widetilde{k}|)\leq\langle \widetilde{A} \rangle$ which is the smallest nonnegative number satisfying
\begin{equation}\label{d2}
\mu[\widetilde{A}]+\rho \langle \widetilde{A} \rangle\geq K.
\end{equation}
Thus, the larger $[\widetilde{A}]$ the more Fourier coefficients are discarded. If $[\widetilde{A}]$  is sufficiently
large the whole $\widetilde{A}$-component is dropped. The upshot is that for the remaining
perturbation $Q$ one has
$$|||P-Q|||_{m-\mu,r-\rho,s}\leq e^{-K}|||P|||_{m,r,s}.$$
Next, each Fourier coefficient of $Q$ is linearized with respect to $z$ at the origin.
Denoting the result of this truncation process by $R$ we obtain
\begin{equation}\label{d3}
|||P-R|||_{m-\mu,r-\rho,\eta s}\leq \Big(e^{-K}+{{\eta^2}\over{1-\eta}}\Big)|||P|||_{m,r,s}
\end{equation}
for $0<\mu<m,0<\rho<r$ and $0<\eta<1$.  Moreover, the estimate
\begin{equation}\label{d0}
|||R|||_{m,r,s}\leq 2|||P|||_{m,r,s}
\end{equation}
obviously holds.

\subsection{Extending the small divisor estimate}
We claim that, if
\begin{equation}\label{e6}
h\leq \min \limits _{\widetilde{A}\in \widetilde{\mathcal{S}}}{{1}\over {\Delta([\widetilde{A}])\langle \widetilde{A}\rangle \Delta(\langle \widetilde{A}\rangle)}}
\end{equation}
with $\langle \widetilde{A}\rangle$ as in the previous subsection \ref{subsec:truncating}, then the estimates
\begin{equation}\label{e000000007}
\begin{aligned}
&{{| \langle k,\omega \rangle+\langle {\widetilde{k}},{\widetilde{\omega}}\rangle|} \geq {1 \over {\Delta\big([[(k,{\widetilde{k}})]]\big)\Delta\big(|k+\widetilde{k}|\big)}}},\\[0.2cm] &\mu[[(k,{\widetilde{k}})]]+\rho (|k|+|\widetilde{k}|)\leq K,\ \ \ (k,{\widetilde{k}}) \in \mathbb{Z}_{\widetilde{\mathcal{S}}}^{{\mathbb{Z}}}\backslash \{\widetilde{k}=0\}
\end{aligned}
\end{equation}
hold uniformly in ${\widetilde{\omega}}$ on the complex neighbourhood $\mathcal{O}_h$ of the set $\mathcal{O}_{\ast}$.

The proof is simple. Given ${\widetilde{\omega}}$ in $\mathcal{O}_h$ there exists an ${\widetilde{\omega}}_0$ in $\mathcal{O}_{\ast}$ such that
$|{\widetilde{\omega}}-{\widetilde{\omega}}_0|<h$.  Given $(k,{\widetilde{k}})$ there exists an $\widetilde{A}$ in $\widetilde{\mathcal{S}}$ containing the support of $(k,{\widetilde{k}})$ such
that $[[(k,{\widetilde{k}})]] = [\widetilde{A}]$. It follows that $|k|+|\widetilde{k}|  \leq\langle \widetilde{A}\rangle$ and hence
\begin{equation*}
| \langle {\widetilde{k}},{\widetilde{\omega}}\rangle-\langle {\widetilde{k}},{\widetilde{\omega}}_0 \rangle| \leq |{\widetilde{k}}|\,
|{\widetilde{\omega}}-{\widetilde{\omega}}_0|\leq \langle \widetilde{A}\rangle h\leq {{1}\over {\Delta([\widetilde{A}])\Delta(\langle \widetilde{A}\rangle)}}
\leq {{1}\over {\Delta([[(k,{\widetilde{k}})]])\Delta(|k|+|\widetilde{k}|)}}
\end{equation*}
by the monotonicity of $\Delta$. The claim follows from the estimate (\ref{d1}) for $\langle {\widetilde{k}},{\widetilde{\omega}}_0 \rangle$.

\subsection{Solving the linearized equation}
The KAM-theorem is proven by the usual Newton-type iteration procedure,
which involves an infinite sequence of coordinate changes and is described in some
detail for example in \cite{Poschel89}. Each coordinate change $\Phi$ is obtained as the time-1-map
$X_F^t|_{t=1}$ of a Hamiltonian vector field $X_F$. Its generating Hamiltonian $F$ as well
as some correction $\widehat{N}$ to the given normal form $N$ are a solution of the linearized
equation
$$\{F,N\}+\widehat{N}=R,$$
which is the subject of this subsection.

The linearized equation $\{F,N\}+\widehat{N}=R$ is broken up into the component equations
$i\big(\langle {k},{\omega}\rangle+\langle {\widetilde{k}},{\widetilde{\omega}}\rangle\big)F_{\widetilde{A}}+\widehat{N}_{\widetilde{A}}=R_{\widetilde{A}}$ with $\mbox{supp}\, (k,\widetilde{k}) \subseteq \widetilde{A}, (k,\widetilde{k})\neq 0$, and solved for $F_{\widetilde{A}}$ and
$\widehat{N}_{\widetilde{A}}$ as described in Section \ref{sec:pro}. Clearly, $\widehat{N}_{\widetilde{A}}=[R_{\widetilde{A}}]$, which is the mean value of $R_{\widetilde{A}}$ over $\mathbb{T}^{\widetilde{A}}$ and $\|\widehat{N}_{\widetilde{A}}\|_{r,s}\leq \|R_{\widetilde{A}}\|_{r,s}$. Hence
\begin{equation}\label{70000}
|||\widehat{N}|||_{m,r,s}\leq |||R|||_{m,r,s}
\end{equation}
by putting pieces together.

The normalized Fourier series expansion of $F_{\widetilde{A}}$ is given by (\ref{100001}). By the extended
small divisor estimate (\ref{e000000007}) and nonresonance condition (\ref{b6}) with $\alpha=\widetilde{\alpha}=2$,
\begin{equation*}
\begin{aligned}
\|F_{\widetilde{A}}\|_{r-\rho,s}&\leq\sum\limits _{\substack{\mbox{supp}\, (k,{\widetilde{k}}) \subseteq \widetilde{A} \\ (k,{\widetilde{k}})\neq 0 }}\Delta\big([[(k,{\widetilde{k}})]]\big)\Delta\big(|k|+|{\widetilde{k}}|\big)|R_{A,k,\widetilde{k}}|_{s}\,e^{(r-\rho)(|k|+|{\widetilde{k}}|)}\\[0.2cm]
&\leq \Delta\big([\widetilde{A}]\big)\Gamma_0(\rho)\|R\|_{r,s},
\end{aligned}
\end{equation*}
where $\Gamma_0(\rho)=\sup \limits _{t\geq 0}\Delta(t)e^{-\rho t}$. Similarly, for the convenience of later estimates,
\begin{equation*}
\begin{aligned}
\sum\limits _{\lambda\in A}\|\partial_{\theta_\lambda}F_{\widetilde{A}}\|_{r-\rho,s}&\leq\sum \limits _{\substack{\mbox{supp}\, (k,{\widetilde{k}}) \subseteq \widetilde{A} \\ (k,{\widetilde{k}})\neq 0 }} \Delta\big([[(k,{\widetilde{k}})]]\big)\,|{k}|\,\Delta\big(|k|+|{\widetilde{k}}|\big)|R_{A,k,\widetilde{k}}|_{s}\,e^{(r-\rho)(|k|+|{\widetilde{k}}|)}\\[0.2cm]
&\leq \Delta\big([\widetilde{A}]\big)\Gamma_1(\rho)\|R\|_{r,s},
\end{aligned}
\end{equation*}
\begin{equation*}
\begin{aligned}
\|\partial_xF_{\widetilde{A}}\|_{r-\rho,s}&\leq\sum \limits _{\substack{\mbox{supp}\, (k,{\widetilde{k}}) \subseteq \widetilde{A} \\ (k,{\widetilde{k}})\neq 0 }} \Delta\big([[(k,{\widetilde{k}})]]\big)\,|{\widetilde{k}}|\,\Delta\big(|k|+|{\widetilde{k}}|\big)|R_{A,k,\widetilde{k}}|_{s}\,e^{(r-\rho)(|k|+|{\widetilde{k}}|)}\\[0.2cm]
&\leq \Delta\big([\widetilde{A}]\big)\Gamma_1(\rho)\|R\|_{r,s},
\end{aligned}
\end{equation*}
where $\Gamma_1(\rho)=\sup \limits _{t\geq 0}(1+t)\Delta(t)e^{-\rho t}$.  Putting the spatial components together,

\begin{equation*}
\begin{aligned}
|||F|||_{m-\mu,r-\rho,s}&\leq\sum\limits _{\widetilde{A}\in\widetilde{\mathcal{S}}}\Delta\big([\widetilde{A}]\big)\Gamma_0(\rho)\|R\|_{r,s}e^{(m-\mu)[\widetilde{A}]}\\[0.2cm]
&\leq \Gamma_0(\mu)\Gamma_0(\rho)|||R|||_{m,r,s}
\end{aligned}
\end{equation*}
and
$$\sum \limits_{\lambda}|||\partial_{\theta_\lambda} F|||_{m-\mu,r-\rho,s},|||\partial_x F|||_{m-\mu,r-\rho,s}\leq \Gamma_0(\mu)\Gamma_1(\rho)|||R|||_{m,r,s}$$
for $0<\mu<m.$

In view of the estimate $\Gamma_0(\rho)\leq \rho \Gamma_1(\rho)$ in  Lemma 6 in \cite{Poschel90}  we may summarize these
estimates by writing
\begin{equation}\label{d5}
\begin{aligned}
\rho^{-1}|||F|||_{m-\mu,r-\rho,s}&\leq\Gamma_\mu\Gamma_\rho|||R|||_{m,r,s},\\[0.2cm]
\sum \limits_{\lambda}|||\partial_{\theta_\lambda} F|||_{m-\mu,r-\rho,s}&\leq\Gamma_\mu\Gamma_\rho|||R|||_{m,r,s},\\[0.2cm]
|||\partial_x F|||_{m-\mu,r-\rho,s}&\leq\Gamma_\mu\Gamma_\rho|||R|||_{m,r,s}
\end{aligned}
\end{equation}
with
\begin{equation}\label{d100009}
\Gamma_\mu=\Gamma_0(\mu)=\sup \limits _{t\geq 0}\Delta(t)e^{-\rho t},\ \ \ \ \Gamma_\rho=\Gamma_1(\rho)=\sup \limits _{t\geq 0}(1+t)\Delta(t)e^{-\rho t}.
\end{equation}

\subsection{The derivatives of F}
On the domain $\mathcal{D}_{r-\rho,s}$  we obtain the estimate
\begin{equation}\label{d6}
\begin{aligned}
|\partial_\theta F|_w=\sum \limits_{\lambda} |\partial_{\theta_\lambda}F|e^{w[\lambda]}&\leq\sum \limits_{\lambda}\sum \limits_{A\ni\lambda}||\partial_{\theta_\lambda}F_{\widetilde{A}}||_{r-\rho,s}e^{w [A]}\\[0.2cm]
&\leq \sum \limits_{\lambda}|||\partial_{\theta_\lambda}F|||_{w,r-\rho,s}
\end{aligned}
\end{equation}
and
\begin{equation}\label{d7}
|\partial_xF|\leq |||\partial_xF|||_{m-\mu,r-\rho,s}.
\end{equation}
Similarly, on the domain $\mathcal{D}_{r-\rho,{s/2}}$ we obtain estimate
\begin{equation}\label{d8}
|\partial_JF|_\infty=0
\end{equation}
and
\begin{equation}\label{d9}
|\partial_zF|\leq |||\partial_zF|||_{m-\mu,r-\rho,{s/2}}\leq {2\over s} |||F|||_{m-\mu,r-\rho,s}.
\end{equation}
Requiring that
\begin{equation}\label{d10}
m-\mu\geq w
\end{equation}
and recalling the estimates (\ref{d01}), (\ref{d0}), (\ref{d5}), (\ref{d6}), (\ref{d7}), (\ref{d8}), (\ref{d9}) we thus have
$${1\over \rho}|\partial_JF|_\infty,{1\over \rho}|\partial_zF|,{2\over s}|\partial_\theta F|_w,{2\over s}|\partial_xF|\leq {2\over s}\Gamma_\mu\Gamma_\rho|||R|||_{m,r,s}\leq 4\,\Gamma_\mu\Gamma_\rho {\varepsilon\over s}$$
uniformly on the domain $\mathcal{D}_{r-\rho,{s/2}}$.

These estimates are expressed more conveniently by means of a weighted phase
space norm. Let
$$|(\theta,x,J,z)|_\mathcal{P}=\max\big(|\theta|_\infty,|x|,|J|_w,|z|),$$
$$W=\text{diag}(\rho^{-1}I_\mathbb{Z},\rho^{-1}I_n,2s^{-1}I_\mathbb{Z},2s^{-1}I_n).$$
Then the above estimates are equivalent to
$$|W X_F|_\mathcal{P}\leq4\,\Gamma_\mu\Gamma_\rho E,\ \ \ \ E={\varepsilon\over s}$$
on $\mathcal{D}_{r-\rho,{s/2}}$.

\subsection{Transforming the coordinates}
The $|W\cdot|_\mathcal{P}$ -distance of the domain
$$\mathcal{D}_\mathfrak{L}=\mathcal{D}_{r-2\rho,{s/4}}\subset\mathcal{D}_L=\mathcal{D}_{r-\rho,{s/2}}$$
to the boundary of $\mathcal{D}_L$ is exactly one half. Hence, if $16\,\Gamma_\mu\Gamma_\rho E\leq 1$, then $|W X_F|_\mathcal{P}$
is less than or equal one fourth on $\mathcal{D}_L$ and consequently
$$X_F^t : \mathcal{D}_\mathfrak{L}\rightarrow \mathcal{D}_L,\ \ \ \ 0\leq t\leq 1.$$
In particular, the time-1-map $\Phi$ is a symplectic map from $\mathcal{D}_{\mathfrak{L}}$ into $\mathcal{D}_L$,  for which the
estimate
\begin{equation}\label{d11}
|W(\Phi-id)|_{\mathcal{P};\mathcal{D}_{\mathfrak{L}}}\leq4\,\Gamma_\mu\Gamma_\rho E
\end{equation}
holds.

In fact, under the present smallness condition on $E$ this statement holds as well
for the larger domain $\mathcal{D}_{r-\kappa\rho,{\kappa s/ 4}}$ instead of $\mathcal{D}_{\mathfrak{L}}$, where $\kappa={3/ 2}$.  The $|W\cdot|_\mathcal{P}$ -distance of its boundary to $\mathcal{D}_{\mathfrak{L}}$ is exactly one fourth. Applying the general Cauchy
inequality of Appendix B in \cite{Poschel90} to the last estimate it follows that in addition,
$$|W(D\Phi-I)W^{-1}|_{\mathcal{P};\mathcal{D}_{\mathfrak{L}}}\leq16\,\Gamma_\mu\Gamma_\rho E,$$
where the norm of derivative is the operator norm induced by $|\cdot|_\mathcal{P}$, see Appendix A in \cite{Poschel90} for its definition.  Finally, if we require
$$4\,\Gamma_\mu\Gamma_\rho E\leq\eta\leq{1\over 2},$$
then
$$X_F^t : \mathcal{D}_\beta=\mathcal{D}_{r-2\rho,{\eta s/ 2}}\rightarrow\mathcal{D}_\eta=\mathcal{D}_{r-\rho,{\eta s}},\ \ \ \ 0\leq t\leq 1$$
by the same arguments as before.

\subsection{Transforming the frequencies}
To put $N_+=N+\widehat{N}$  into normal form, the frequency parameters are transformed by setting ${\widetilde{\omega}}_+={\widetilde{\omega}}+ v({\widetilde{\omega}})$.
Proceeding just as in (\ref{70000}) the estimate for $\widehat{N}$
implies that $|v|_{h/2}=|\partial_z\widehat{N}|_{h/ 2}\leq 4E$.  Referring to Lemma 11 in \cite{Poschel90} or  Lemma A.3 in \cite{Poschel01} it follows that for
\begin{equation}\label{e15}
E\leq {h\over 16}
\end{equation}
and the map $id+v$  has a real analytic inverse
$$\varphi : \ \ {\mathcal{O}_{h/ 4}}\rightarrow {\mathcal{O}_{h/ 2}},\ \ {\widetilde{\omega}}_+\mapsto{\widetilde{\omega}}$$
with the estimate
\begin{equation}\label{e1006}
|\varphi-id|,\ \ \ {h\over 4}\Big|D \varphi-Id\Big|\leq 4E
\end{equation}
uniformly on $\mathcal{O}_{h/4}$.

\subsection{Estimating the new error term}
The new error term is
$$P_+=\int_0^1\{R_t,F\}\circ X_F^t\, dt+(P-R)\circ X_F^1,$$
where $R_t=(1-t)\widehat{N}+tR$. By Lemma 10 in \cite{Poschel90} and estimate (\ref{d5}),
$$|||G\circ X_F^t|||_{m-\mu,r-2\rho,{{\eta s}/ 2}}\leq 2|||G|||_{m-\mu,r-\rho,{\eta s}},\ \ \ \ 0\leq t\leq 1,$$
provided that
\begin{equation}\label{e16}
4\,C_0\,\Gamma_\mu\Gamma_\rho E\leq\eta\leq{1\over2},
\end{equation}
where $C_0=8$  is a constant. Hence, with this assumption,
$$|||P_+|||_{m-\mu,r-2\rho,{{\eta s}/ 2}}\leq 2|||\{R_t,F\}|||_{m-\mu,r-\rho,{\eta s}}+2|||(P-R)|||_{m-\mu,r-\rho,{\eta s}}.$$
Obviously, $|||R_t|||_{m,r,s}\leq 2\varepsilon$ for $0\leq t\leq 1$  by the estimates for
$\widehat{N}$ and $F$, and therefore by (\ref{d01}), (\ref{d0}), (\ref{d5}), we get
\begin{equation*}
|||\partial_xR_{t}|||_{m-\mu,r-\rho,{\eta s}}\, |||\partial_zF|||_{m-\mu,r-\rho,{\eta s}}\leq {{2\varepsilon}\over \rho} {2\over s} |||F|||_{m-\mu,r-\rho,s}
\leq 8\,\Gamma_\mu\Gamma_\rho E \varepsilon,
\end{equation*}
\begin{equation*}
|||\partial_zR_{t}|||_{m-\mu,r-\rho,{\eta s}}\, |||\partial_xF|||_{m-\mu,r-\rho,{\eta s}}\leq  {2\over s} 2\varepsilon|||\partial_xF|||_{m-\mu,r-\rho,s}
\leq 8\,\Gamma_\mu\Gamma_\rho E \varepsilon.
\end{equation*}
Hence
\begin{equation*}
\begin{aligned}
|||\{R_t,F\}|||_{m-\mu,r-\rho,{\eta s}}&\leq |||\partial_xR_{t}|||_{m-\mu,r-\rho,{\eta s}}\, |||\partial_zF|||_{m-\mu,r-\rho,{\eta s}}\\[0.2cm]
&+|||\partial_zR_{t}|||_{m-\mu,r-\rho,{\eta s}}\, |||\partial_xF|||_{m-\mu,r-\rho,{\eta s}}\\[0.2cm]
&\leq 16\,\Gamma_\mu\Gamma_\rho E \varepsilon
\end{aligned}
\end{equation*}
in view of (\ref{d8}) and $R_t$ independent of $J$. Combined with (\ref{d01}), (\ref{d3}) we altogether
obtain
\begin{equation}
|||P_+|||_{m-\mu,r-2\rho,{{\eta s}/ 2}}\leq 32\,\Gamma_\mu\Gamma_\rho E \varepsilon+2 e^{-K}\varepsilon +4\eta^ 2 \varepsilon
\end{equation}
for the new error term.

\section{Iteration and Convergence}\label{sec:ite}
\subsection{The iterative construction}
To iterate the KAM step infinitely often we now choose sequences for the pertinent parameters. Let $a=13,b=4,c=6,d=8,e=22$ and $\kappa={3/ 2}$. The choice of these integer
constants will be motivated later in the course of the proof of the iterative lemma.

Given $0<\mu\leq m-w$ and $0<\rho<{r/ 2}$  there exist sequences $\mu_0\geq\mu_1\geq\cdots>0$ and $\rho_0\geq\rho_1\geq\cdots>0$ such that
\begin{equation}\label{f10001}
\Psi_0(\mu)\Psi_1(\rho)=\prod_{\nu=0}^\infty\Gamma_{\mu_\nu}^{\kappa_\nu}\Gamma_{\rho_\nu}^{\kappa_\nu}
\end{equation}
with
$$\sum_{\nu=0}^\infty \mu_\nu=\mu,\ \ \ \ \sum_{\nu=0}^\infty \rho_\nu=\rho,\ \ \ \ \kappa_\nu={{\kappa-1}\over{\kappa^{\nu+1}}}$$
where $\Gamma_{\mu_\nu}=\Gamma_{0}(\mu_\nu)$ and $\Gamma_{\rho_\nu}=\Gamma_{1}(\rho_\nu)$, $\Gamma_{0},\Gamma_{1}$ is defined by (\ref{d100009}).  Fix such sequences, and for $j\geq 0$ set
$$\Gamma_j=2^{j+a}\Gamma_{\mu_j}\Gamma_{\rho_j},\ \ \ \Theta_j=\prod_{\nu=0}^{j-1}\Gamma_{\nu}^{\kappa_\nu}, \ \ \ E_j=(\Theta_j E_0)^{\kappa^j},$$
where $\Theta_0=1$.  Furthermore, set
$$m_j=m-\sum_{\nu=0}^{j-1}\mu_j,\ \ \ \ r_j=r-2\sum_{\nu=0}^{j-1}\rho_\nu,$$
\begin{equation}\label{f101}
s_j=s\prod_{\nu=0}^{j-1}{{\eta_\nu}\over 2},\ \ \ \ h_j=2^{j+c}E_j,
\end{equation}
where $\eta_j^2=4^{-b}\Gamma_jE_j$. Then $m_j\downarrow m-w, r\downarrow r-2\rho$ and $s_j\downarrow 0$, $h_j\downarrow 0$. These
sequences define the complex domains
$$\mathcal{D}_j=\mathcal{D}_{r_j,s_j},\ \ \ \ \mathcal{O}_j=\mathcal{O}_{h_j}.$$
Finally, we introduce an extended phase space norm,
$$|(\theta,x,J,z,{\widetilde{\omega}})|_\mathcal{\bar{P}}=\max\big(|\theta|_\infty,|x|,|J|_w,|z|,|{\widetilde{\omega}}|),$$
and the corresponding weight matrices,
$$\bar{W}_j=\text{diag}(\rho^{-1}_jI_\mathbb{Z},\rho^{-1}_jI_n,2s^{-1}_jI_\mathbb{Z},2s^{-1}_jI_n,h^{-1}_jI_n).$$

Then we can state the Iterative Lemma.

\begin{lemma}[Iterative Lemma]\label{lem6.1}
Suppose that
\begin{equation}\label{f102}
s^{-1}|||H-N|||_{m,r,s,h}\leq{{\widetilde{\alpha} \varepsilon_*}\over \Psi_0(\mu)\Psi_1(\rho)}\leq {h\over{2^c}},
\end{equation}
where $\widetilde{\alpha}=2$ and $\varepsilon_*=2^{-e}$. Then for each $j\geq 0$ there exists a normal form $N_j$
and a real analytic transformation
$$\mathcal{F}^j = \mathcal{F}_0\circ\mathcal{F}_1\circ\cdots\circ\mathcal{F}_{j-1} :\ \ \ \mathcal{D}_j\times\mathcal{O}_j\rightarrow \mathcal{D}_0\times\mathcal{O}_0$$
of the form described in Section \ref{sec:pro}, which is symplectic for each ${\widetilde{\omega}}$, such that $H\circ\mathcal{F}^j=N_j+P_j$ with
\begin{equation}\label{f103}
s_j^{-1}|||P_j|||_{m_j,r_j,s_j,h_j}\leq E_j.
\end{equation}
Moreover,
\begin{equation}\label{f1004}
|\bar{W}_0(\mathcal{F}^{j+1}-\mathcal{F}^j)|_\mathcal{\bar{P}}\leq4\,\max\Big(2^{1-a-j}\Gamma_j E_j,2 {{E_j}/{h_j}}\Big)
\end{equation}
on $\mathcal{D}_{j+1}\times\mathcal{O}_{j+1}$.
\end{lemma}

Before giving the proof of Iterative Lemma \ref{lem6.1} we collect some useful facts. The $\kappa_\nu$
satisfy the identities
$$\sum_{\nu=0}^\infty \kappa_\nu=1,\ \ \ \ \sum_{\nu=0}^\infty \nu\kappa_\nu={1\over{\kappa-1}}.$$
This and the monotonicity of the $\Gamma$ -function imply that
$$\Gamma_j=\prod_{\nu=j}^\infty\Gamma_{j}^{\kappa_\nu \kappa^j}\leq\Bigg(\prod_{\nu=j}^\infty\Gamma_{j}^{\kappa_\nu }\Bigg)^{\kappa^j}.$$
Together with the definition of $E_j$ and (\ref{f10001}) we obtain the estimate
\begin{equation}\label{f2}
\Gamma_jE_j\leq \Bigg(\prod_{\nu=0}^\infty\Gamma_{\nu}^{\kappa_\nu }E_0\Bigg)^{\kappa^j}=\big(2^{2+a}\Psi_0\Psi_1E_0)\big)^{\kappa^j}.
\end{equation}
Moreover,
\begin{equation}\label{f3}
\Gamma_j^{\kappa-1}E_j^\kappa=E_{j+1}
\end{equation}
by a straightforward calculation.

\noindent\textbf{Proof of Iterative Lemma \ref{lem6.1}:}  Lemma \ref{lem6.1}  is proven by induction. Choosing $\mathcal{F}^0= id$ and
$$E_0={{\widetilde{\alpha} \varepsilon_*}\over \Psi_0(\mu)\Psi_1(\rho)},$$
there is nothing to prove for $j=0$. Just observe that $h_0\leq h$ by the very definition
of $h_0$ and $E_0$.

Let $j\geq 0$.  To apply the KAM-step to $H_j=H\circ\mathcal{F}^j$ and $N_j$ we need to verify
its assumptions (\ref{e6}), (\ref{d10}), (\ref{e15}) and (\ref{e16}). Clearly, $m_j-\mu_j\geq w$  by construction,
and $E_j\leq {{h_j}/ 16}$  in view of the definition of $h_j$ and $c\geq 4$,  so the second and third
requirements are met. Taking squares, the fourth requirement is equivalent to
$$4^{2-a-j}C_0^2\Gamma_j^2E_j^2\leq 4^{-b}\Gamma_jE_j\leq{1\over 4}.$$
This holds for all $j\geq0$, since $C_0=8$,
\begin{equation}\label{f4}
\Gamma_jE_j\leq2^{3+a-e}
\end{equation}
by (\ref{f101}), (\ref{f102}), (\ref{f2}) and $a \geq b + 2, b \geq 0, e \geq a + 9$.

As to the first requirement, define $K_j$ by $e^{-K_j}=2^{-d}\,\Gamma_jE_j$ and subsequently
$\langle \cdot\rangle$ as in (\ref{d2}). For arbitrary $\widetilde{A}$ in $\widetilde{\mathcal{S}}$ with $\langle\widetilde{ A}\rangle>0$ we then have
\begin{equation}
\begin{aligned}
{1\over{\langle \widetilde{A}\rangle \Delta(\langle \widetilde{A}\rangle) \Delta([\widetilde{A}])}}&= {{e^{-\rho_j\langle \widetilde{A}\rangle}e^{-\mu_j[ \widetilde{A}]}}\over{\langle \widetilde{A}\rangle \Delta(\langle \widetilde{A}\rangle) e^{-\rho_j\langle \widetilde{A}\rangle} \Delta([\widetilde{A}])e^{-\mu_j[ \widetilde{A}]}}}\\[0.2cm]
&\geq {{e^{-K_j}}\over{\Gamma_{\mu_j}\Gamma_{\rho_j}}}={{2^{-d}\Gamma_jE_j}\over {2^{-j-2}\Gamma_j}}=2^{j+a-d}E_j\geq h_j,
\end{aligned}
\end{equation}
since $a \geq c + d$. This estimate holds even more when $\langle \widetilde{A}\rangle=0$. Hence, also
requirement (\ref{e6}) is satisfied.

The KAM-construction now provides a normal form $N_{j+1}$, a coordinate transformation $\Phi_j$ and a parameter transformation $\varphi_j$. By the definition of $r_j$ and $s_j$,
$\Phi_j$ maps $\mathcal{D}_{j+1}$ into $\mathcal{D}_{j}$, while $\varphi_j$  maps $\mathcal{O}_{j+1}$ into $\mathcal{O}_{j}$, since
$${{h_{j+1}\over {h_j}}}={{2E_{j+1}}\over E_j}=2(\Gamma_jE_j)^{\kappa-1}\leq2^{1+(3+a-e)/2}\leq{1\over 4}$$
in view of (\ref{f4}) and $\kappa={3/2},e\geq a+9$. Setting
$$\mathcal{F}^{j+1}=\mathcal{F}^j\circ\mathcal{F}_j,\ \ \ \ {F}_j=\Phi_j\circ\varphi_j,$$
we obtain a transformation $\mathcal{F}^{j+1}$ from $\mathcal{D}_{j+1}\times\mathcal{O}_{j+1}$ into $\mathcal{D}_0\times\mathcal{O}_0$. For the new
error term
$$P_{j+1}=H\circ\mathcal{F}^{j+1}-N_{j+1}=H_j\circ{F}_j-N_{j+1},$$
we obtain
\begin{equation*}
\begin{aligned}
|||P_{j+1}|||_{j+1}&\leq 32\,\Gamma_{\mu_j}\Gamma_{\rho_j} E_j \varepsilon_j+2 e^{-K_j}\varepsilon_j +4\eta_j^ 2 \varepsilon_j\\[0.2cm]
&\leq (2^{5-a}+2^{1-d}+2^{2-2b})\Gamma_jE_j\varepsilon_j.
\end{aligned}
\end{equation*}
Dividing by $s_{j+1}=\eta_js_j/2$ this yield
\begin{equation*}
\begin{aligned}
s^{-1}_{j+1}|||P_{j+1}|||_{j+1}&\leq 2^{1+b}(2^{5-a}+2^{1-d}+2^{2-2b})\Gamma_j^{\kappa-1}E_j^\kappa\\[0.2cm]
&= (2^{6-a+b}+2^{2+b-d}+2^{3-b})E_{j+1}\\[0.2cm]
&\leq E_{j+1},
\end{aligned}
\end{equation*}
since $\eta_j^2=4^{-b}\Gamma_jE_j, E_j=\varepsilon_j/s_j, a \geq b + 8, b \geq 4, d \geq b + 4, \kappa={3/ 2}$  and (\ref{f3}).

To prove the first of the estimates, write
\begin{equation}\label{f9}
\begin{aligned}
|\bar{W}_0(\mathcal{F}^{j+1}-\mathcal{F}^j)|_{j+1}&=|\bar{W}_0(\mathcal{F}^{j}\circ\mathcal{F}_j-\mathcal{F}^j)|_{j+1}\\[0.2cm]
&\leq |\bar{W}_0\bar{D}\mathcal{F}^{j}\bar{W}_j^{-1}|_{j}\, |\bar{W}_j(\mathcal{F}_j-id)|_{j+1},
\end{aligned}
\end{equation}
where $|\cdot|_j=|\cdot|_{\mathcal{\bar{P}},\mathcal{D}_j\times\mathcal{O}_j}$, where $\bar{D}$  denotes differentiation with respect to $(x,\theta,z,J,{\widetilde{\omega}})$. By (\ref{d11}), (\ref{e1006}) and the definition of $\Gamma_j$,
\begin{equation}\label{f10}
\begin{aligned}
|\bar{W}_j(\mathcal{F}_j-id)|_{j+1}&\leq \max\big(|W_j(\Phi_j-id)|_\mathcal{P},h^{-1}_j|\varphi-id|\big)\\[0.2cm]
&\leq \max\big(2^{2-a-j}\Gamma_j E_j,4{{E_j}/ {h_j}}\big).
\end{aligned}
\end{equation}
It remains to show that the first factor is bounded by 2. By the inductive construction,
$\mathcal{F}^j = \mathcal{F}_0\circ\mathcal{F}_1\circ\cdots\circ\mathcal{F}_{j-1}$, and
\begin{equation*}
\begin{aligned}
|\bar{W}_\nu \bar{D}\mathcal{F}_\nu \bar{W}_\nu^{-1}|_{\nu+1}&\leq \max\big(|{W}_\nu \bar{D}\Phi_\nu {W}_\nu^{-1}|_\mathcal{P}+h_\nu|W_\nu\partial_{{\widetilde{\omega}}}\Phi_\nu|_\mathcal{P},|\partial_{{\widetilde{\omega}}}\varphi_\nu|\big)\\[0.2cm]
&\leq \max\big(1+2^{5-a-\nu}\Gamma_\nu E_\nu,1+16{{E_\nu}/ {h_\nu}}\big)\\[0.2cm]
&\leq 1+2^{4-c-\nu}.
\end{aligned}
\end{equation*}
By (\ref{d11}), (\ref{e1006}), (\ref{f4}).  Since the weights of $\bar{W}_\nu^{-1}$ do not decrease as $\nu$ decreases, and
since $c\geq 6$, we obtain
\begin{equation}\label{f11}
\begin{aligned}
|\bar{W}_0 \bar{D}\mathcal{F}^j \bar{W}_j^{-1}|_{j}\leq\prod_{\nu=0}^{j-1}|\bar{W}_\nu \bar{D}\mathcal{F}_\nu \bar{W}_{\nu+1}^{-1}|_{\nu+1}\leq \prod_{\nu=0}^{\infty}(1+2^{4-c-\nu})\leq 2.
\end{aligned}
\end{equation}
By (\ref{f9}), (\ref{f10}), (\ref{f11}), the conclusion (\ref{f1004}) holds. This completes
the proof of the Iterative Lemma \ref{lem6.1}.\qed

\subsection{Convergence}
By the estimates of the iterative lemma the $\mathcal{F}^j$ converge uniformly on
$$\bigcap\limits _{j\geq 0}\mathcal{D}_{j}\times\mathcal{O}_{j}=\mathcal{D}_*\times\mathcal{O}_{\ast},\ \ \  \ \mathcal{D}_*=\mathcal{D}_{r-2\rho,0}$$
to mappings $\mathcal{F}_*$  that are real analytic in $x,\theta$ and uniformly continuous in ${\widetilde{\omega}}$.
Moreover,
$$|\bar{W}_0(\mathcal{F}_*-id)|_{\bar{\mathcal{P}}}\leq {1\over 2}$$
on $\mathcal{D}_*\times\mathcal{O}_{\ast}$ by the usual telescoping argument.

But by construction, the $\mathcal{F}^j$ are affine linear in each fiber over $\mathbb{T}^\mathbb{Z}\times\mathbb{T}^n\times\mathcal{O}_{\ast}$. Therefore they indeed converge uniformly on any domain $\mathcal{D}_{r-2\rho,\sigma}\times\mathcal{O}_{\ast}$ with $\sigma>0$ to a map $\mathcal{F}_*$ that is real analytic and symplectic for each ${\widetilde{\omega}}$. In particular,
$$\mathcal{F}_* : \mathcal{D}_{r-2\rho,{s/ 2}}\times\mathcal{O}_{\ast}\rightarrow \mathcal{D}_{r,s}\times\mathcal{O}_h$$
by piecing together the above estimates.

Going to the limit in (\ref{f103}) and using Cauchy's inequality we finally obtain
$$H\circ\mathcal{F}_*=e_*+\langle {\omega},J\rangle+\langle {\widetilde{\omega}},z\rangle+\cdots.$$
This completes the proof of Theorem  \ref{thm3.2}.

\subsection{Estimates}
The scheme so far provides only a very crude estimate of $\mathcal{F}_*$ since the actual size
of the perturbation is not taken into account in the estimates of the iterative lemma.
But nothing changes when all inequalities are scaled down by the factor ${\varepsilon/ E}\leq 1$,
where
$$\varepsilon=s^{-1}|||H-N|||_{m,r,s,h}\leq E={{\alpha \varepsilon_*}\over \Psi_\mu\Psi_\rho}.$$
It follows that
$$|\bar{W}_0(\mathcal{F}_{*}-id|_\mathcal{\bar{P}}\leq{\varepsilon\over E}$$
uniformly on $\mathcal{D}_{r-2\rho,{s/ 2}}\times\mathcal{O}_\ast$.

\section{The measure estimate}\label{sec:mea}
In this section the measure estimate of the frequency ${\widetilde{\omega}}$ satisfying inequalities  (\ref{c000000004}) will be given. Firstly, we give some useful  lemmas.
\begin{lemma}[Lemma\ 2\ in {\cite{Poschel90}}]\label{lem4.1}
For every given approximation function $\Theta$, there exists an approximation
function $\Delta$ such that
$$\sum \limits_{\widetilde{A}\in\widetilde{\mathcal{S}},\mbox{card} (\widetilde{A})=i} {1\over {\Delta([\widetilde{A}])}}\leq {{2N_0}\over {\Theta(t_i)}},\ \ \ i\geq 1,$$
where $N_0, t_i$ are given in Lemma \ref{lem2.11}.
\end{lemma}

\begin{lemma}[Lemma\ 4\ in {\cite{Poschel90}}]\label{lem4.3}
There is an approximation function $\Delta$ such that
$$\sum \limits_{\ell\in\mathbb{Z}^i\backslash \{0\}} {1\over {\Delta(|\ell|)}}\leq {M}^{i\log\log i}$$
for all sufficiently large $i$ with some constant ${M}$, where $\ell=(\ell_1,\ell_2,\cdots,\ell_i)$ and $|\ell|=|\ell_1|+|\ell_2|+\cdots+|\ell_i|$.
\end{lemma}

\begin{remark}\label{rem7.3}
Of course, Lemma \ref{lem4.3} also gives a bound for all small $i$, since the left hand side is monotonically increasing with $i$.
\end{remark}

\begin{remark}
From the proofs of Lemma \ref{lem4.1} and  Lemma \ref{lem4.3}, we know that there is the same approximation function $\Delta$ such that Lemma \ref{lem4.1} and  Lemma \ref{lem4.3} hold simultaneously. The detail proofs of Lemma \ref{lem4.1} and  Lemma \ref{lem4.3} can be found in {\cite{Poschel90}}.
\end{remark}

\begin{theorem}\label{thm2.12}
There is an approximation function $\Delta$ such that for suitable $\alpha$ ,\ the set $\mathcal{O}_{\alpha}$ of ${\widetilde{\omega}}$ satisfying (\ref{c000000004}) has positive measure.
\end{theorem}
\Proof Choose the frequency  ${\omega}=(\cdots,{\omega}_\lambda,\cdots)$  satisfying the nonresonance condition (\ref{b6}). For any bounded $\mathcal{O} \in\mathbb{R}^n$, let $\mathcal{O}_{\alpha}\subseteq \mathcal{O}$ denote the set of all ${\widetilde{\omega}}$ satisfying (\ref{c000000004}) with fixed $\omega$. Then complement of the open dense set $\mathcal{R}_{\alpha}$, where
\begin{eqnarray*}
\mathcal{R}_{\alpha}&=& \bigcup \limits_{(k,{\widetilde{k}}) \in \mathbb{Z}_{\widetilde{\mathcal{S}}}^{{\mathbb{Z}}}\backslash \{\widetilde{k}=0\}} \mathcal{R}_{\alpha}^{{k},{\widetilde{k}}}\\[0.2cm]
&=&\bigcup \limits_{ (k,{\widetilde{k}}) \in \mathbb{Z}_{\widetilde{\mathcal{S}}}^{{\mathbb{Z}}}\backslash \{\widetilde{k}=0\} } \Bigg\{{\widetilde{\omega}} \in \mathbb{R}^n : \Big|\langle {k},{\omega} \rangle+\langle {\widetilde{k}},{\widetilde{\omega}}\rangle\Big|<{\alpha \over {\Delta\big([[(k,\widetilde{k})]]\big)\Delta\big(|{k}|+|{\widetilde{k}}|\big)}}\Bigg\}.
\end{eqnarray*}

Now we estimate the measure of the set $\mathcal{R}_{\alpha}^{{k},{\widetilde{k}}}$. Since  $(k,{\widetilde{k}}) \in \mathbb{Z}_{\widetilde{\mathcal{S}}}^{{\mathbb{Z}}}\backslash \{\widetilde{k}=0\}$, then $\widetilde{k}\neq0$,  set $|{\widetilde{k}}_{\max}|=\max\limits_{1\leq \imath\leq n}|{\widetilde{k}}_\imath|\neq 0$, then there exists some  $1\leq \jmath\leq n$ such that $|{\widetilde{k}}_\jmath|=|{\widetilde{k}}_{\max}|$. Therefore, we have

\begin{eqnarray*}
\mathcal{R}_{\alpha}^{{k},{\widetilde{k}}}&=&\Bigg\{{\widetilde{\omega}} \in \mathbb{R}^n : \Big|\langle {k},{\omega} \rangle+\langle {\widetilde{k}},{\widetilde{\omega}}\rangle \Big|<{\alpha \over {\Delta\big([[(k,\widetilde{k})]]\big)\Delta\big(|{k}|+|{\widetilde{k}}|\big)}}\Bigg\}\\[0.2cm]
&=& \Bigg\{{\widetilde{\omega}} \in \mathbb{R}^n :\Big|{\widetilde{k}}_{\max}{\widetilde{\omega}}_\jmath+\sum \limits _{\imath\not =\jmath} {\widetilde{k}}_i{\widetilde{\omega}}_i+\langle {k},{\omega}\rangle\Big|<{\alpha \over {\Delta\big([[(k,\widetilde{k})]]\big)\Delta\big(|{k}|+|{\widetilde{k}}|\big)}}\Bigg\}\\[0.2cm]
&=& \Big\{{\widetilde{\omega}} \in \mathbb{R}^n : |{\widetilde{\omega}}_\jmath+b_{{k},{\widetilde{k}}}|<\delta_{{k},{\widetilde{k}}}\Big\}\\[0.2cm]
&=& \Big\{{\widetilde{\omega}} \in \mathbb{R}^n :-b_{{k},{\widetilde{k}}}-\delta_{{k},{\widetilde{k}}}<{\widetilde{\omega}}_\jmath<-b_{{k},{\widetilde{k}}}+\delta_{{k},{\widetilde{k}}}\Big\},
\end{eqnarray*}
where $b_{{k},{\widetilde{k}}}={1\over {{\widetilde{k}}_{\max}}}\Bigg(\sum \limits _{\imath\not =\jmath} {\widetilde{k}}_\imath{\widetilde{\omega}}_\imath+\langle {k},{\omega}\rangle\Bigg)$ and
$
\delta_{{k},{\widetilde{k}}} = {\alpha \over {\Delta\big([[(k,\widetilde{k})]]\big)\Delta\big(|{k}|+|{\widetilde{k}}|\big)}}\,{1\over {|{\widetilde{k}}_{\max}|}}.
$
Obviously, for any bounded domain $\mathcal{O}\subset\mathbb{R}^n$,  we have the Lebesgue measure estimate
$$
\mbox{meas}\big(\mathcal{R}_{\alpha}^{{k},{\widetilde{k}}}\cap\mathcal{O}\big)\le C_0\delta_{{k},{\widetilde{k}}} = {C_0\alpha \over {\Delta\big([[(k,\widetilde{k})]]\big)\Delta\big(|{k}|+|{\widetilde{k}}|\big)}}\,{1\over {|{\widetilde{k}}_{\max}|}}
$$
with some positive constant $C_0$.

Since $|{\widetilde{k}}_{\max}|\neq 0, {\widetilde{k}}_{\max}\in \mathbb{Z}$, which means $|{\widetilde{k}}_{\max}|\geq 1$, then we have the following measure estimate
$$\alpha^{-1}\mbox{meas}\big(\mathcal{R}_{\alpha}^{{k},{\widetilde{k}}}\cap\mathcal{O}\big)\leq{C_0 \over {\Delta\big([[(k,\widetilde{k})]]\big)\Delta\big(|{k}|+|{\widetilde{k}}|\big)}}.$$

Next we estimate the measure of the set $\mathcal{R}_{\alpha}$. From the definition of $\mathbb{Z}_{\widetilde{\mathcal{S}}}^{{\mathbb{Z}}}$, there exists a nonempty  set $\widetilde{A}\in \mathcal{S}$ such that $\mbox{supp}\, (k,\widetilde{k})\subseteq \widetilde{A}$, we get
\begin{align*}
\alpha^{-1}\mbox{meas}(\mathcal{R}_{\alpha}\cap\mathcal{O}) &\leq  \sum \limits _{\widetilde{A}\in \widetilde{\mathcal{S}}}\sum \limits _{\substack{\mbox{supp}\,(k,\widetilde{k})\subseteq \widetilde{A}\\ \widetilde{k}\neq 0}}\alpha^{-1}\mbox{meas}\big(\mathcal{R}_{\alpha}^{{k},{\widetilde{k}}}\cap\mathcal{O}\big)\nonumber\\[0.2cm]
&\leq C_0\sum \limits _{\widetilde{A}\in \widetilde{\mathcal{S}}}\sum \limits _{\substack{\mbox{supp}\,(k,\widetilde{k})\subseteq \widetilde{A}\\ \widetilde{k}\neq 0}}{1 \over {\Delta\big([[(k,\widetilde{k})]]\big)\Delta\big(|{k}|+|{\widetilde{k}}|\big)}}\nonumber\\[0.2cm]
&\leq  C_0 \sum\limits_{\widetilde{A}\in\mathcal{S}}\Bigg({{1}\over \Delta([\widetilde{A}])}\sum\limits_{\substack{\mbox{supp}\,(k,\widetilde{k})\subseteq \widetilde{A}\\ \widetilde{k}\neq 0}}{1\over {\Delta\big(|{k}|+|{\widetilde{k}}|\big)}}\Bigg)\nonumber\\[0.2cm]
&\leq  C_0 \sum\limits_{\widetilde{A}\in\mathcal{S}}\Bigg({{1}\over \Delta([\widetilde{A}])}\sum\limits_{\substack{\mbox{supp}\,(k,\widetilde{k})\subseteq \widetilde{A}\\ (k,\widetilde{k})\neq 0}}{1\over {\Delta\big(|{k}|+|{\widetilde{k}}|\big)}}\Bigg)\nonumber\\[0.2cm]
&\leq  C_0\sum\limits_{i=1}^{+\infty}\Bigg(\Bigg(\sum\limits_{\widetilde{A}\in\mathcal{S},\mbox{card} (\widetilde{A})=i}{1\over \Delta([\widetilde{A}])}\Bigg)\sum\limits_{\ell\in \mathbb{Z}^{i}\backslash\{0\}}{1\over {\Delta(|\ell|)}}\Bigg).
\end{align*}
Thus the sum is broken up with respect to the cardinality and the weight of the spatial
components of $\mathcal{S}$. Each of these factors is now studied separately.

By applying Lemma \ref{lem4.1}, Lemma \ref{lem4.3} and Remark \ref{rem7.3}, we arrive at
\begin{eqnarray*}
\alpha^{-1}\mbox{meas}(\mathcal{R}_{\alpha}\cap\mathcal{O}) 
\leq C+C\sum\limits_{i=i_0}^{+\infty}{{{M}^{i\log\log i}}\over{\Theta(t_i)}}
\end{eqnarray*}
with some constant $C$ and $i_0$ so large that $t_i\geq i\log^{\varrho-1} i$ for $i\geq i_0, \varrho>2$ by hypotheses. Here we are still free to choose a suitable approximation function $\Theta$, and choose
$$\Theta(t)=\exp\Bigg({t\over{\log\,t\log^{\varrho-1}\log\,t}}\Bigg),\ \ \ t>e, \varrho>2,$$
the infinite sum does converge. Thus there is an approximation function $\Delta$ such that
$$\alpha^{-1}\mbox{meas}(\mathcal{R}_{\alpha}\cap\mathcal{O})<+\infty.$$
Hence,
$$\mbox{meas}(\mathcal{R}_{\alpha}\cap\mathcal{O})\leq O(\alpha)$$
and
$$\mathcal{O}_{\alpha}\to \mathcal{O}\ \ \ \ \ \mbox{as} \ \ \ \ \ \alpha\to 0.$$
This completes the proof of Theorem \ref{thm2.12}.\qed

\section{Application}\label{sec:app}
In this section we will apply Theorem \ref{thm3.2}  to the differential equation with superquadratic potentials depending almost periodically on time
\begin{equation}\label{h1}
\ddot{x}+x^{2l+1}=\sum_{j=0}^{2l} p_j(t)x^j,
\end{equation}
where $p_0,p_1,\cdots,p_{2l}$ are real analytic almost periodic functions with the frequency ${{\omega}}=(\cdots,{{\omega}}_\lambda,\cdots)$ and admit a  spatial series expansion similar to (\ref{b10004}).

\subsection{Rescaling}
We first rescale the time variable $t$ and the space variable $x$ to get a slow system. Let $u=\varepsilon x, \tau=\varepsilon^{-l}t$. Then
equation (\ref{h1}) becomes
\begin{equation}\label{h2}
u''+u^{2l+1}=\varepsilon\sum_{j=0}^{2l} \varepsilon^{2l-j} p_j(\tau)u^j,
\end{equation}
where  $''$ stands for ${d^2}\over {d\tau^2}$, $p_0,p_1,\cdots,p_{2l}$ are real analytic almost periodic functions in $\tau$ with the frequency $\widehat{\omega}=\varepsilon^l {\omega}$. Without causing confusion, in the following we still use $t$ instead of $\tau$. Equation (\ref{h2}) is equivalent to the following Hamiltonian system
\begin{equation}\label{h3}
\left\{\begin{array}{ll}
u'=v, \\[0.1cm]
v'=-u^{2l+1}+\varepsilon\dsum _{j=0}^{2l}\varepsilon^{2l-j} p_j(t)u^j,
\end{array}\right.
\end{equation}
and the corresponding Hamiltonian function is
\begin{equation}\label{h4}
h(u,v,t)={1\over 2}v^2+{1\over{2l+2}}u^{2l+2}-\varepsilon\sum _{j=0}^{2l} {\varepsilon^{2l-j}\over{j+1}} p_j(t)u^{j+1}.
\end{equation}
It is obvious that (\ref{h4}) is a perturbation of the integrable Hamiltonian
\begin{equation}\label{h5}
h_0(u,v,t)={1\over 2}v^2+{1\over{2l+2}}u^{2l+2},
\end{equation}
for $\varepsilon>0$ small. Our aim is to construct, for every sufficiently small $\varepsilon>0$, invariant cylinders
tending to the infinity for (\ref{h4}) close to $\{h_0(u,v,t)=C\}\times\mathbb{R}$ in the extended phase space, which
prohibit any solution from going to the infinity. For this purpose, we will introduce the action-angle variables first.

\subsection{Action and angle variables}
We consider the following integrable Hamiltonian system
\begin{equation}\label{h6}
\left\{\begin{array}{ll}
u'=v, \\[0.1cm]
v'=-u^{2l+1},
\end{array}\right.
\end{equation}
with Hamiltonian function (\ref{h5}). Suppose $(C(t),S(t))$ is the solution of (\ref{h6}) satisfying the initial condition $(C(0),S(0))=(1,0)$. Let $T_*>0$ be its minimal period, which is a constant. Then these analytic functions $C(t),S(t)$ satisfy\\[0.2cm]
(i) $C(t+T_*)=C(t), S(t+T_*)=S(t)\ \text{and}\ C(0)=1, S(0)=0;$\\[0.2cm]
(ii) $\dot{C}(t)=S(t), \dot{S}(t)=-C^{2l+1}(t);$\\[0.2cm]
(iii) $(l+1)S^2(t)+C^{2l+2}(t)=1;$\\[0.2cm]
(iiii) $C(-t)=C(t), S(-t)=-S(t).$

The action and angle variables are now defined by the map $\Psi :\mathbb{R}^+\times \mathbb{T}\rightarrow \mathbb{R}^2\setminus \{0\}$ via $(\varrho,\phi)=\Psi(u,v)$ is given by the formula
\begin{equation*}
\Psi :\begin{array}{ll}
u=c^{1\over {l+2}}_1{\varrho^{1\over {l+2}}}C\big( {{{T_*}}\phi\over {2\pi}}\big), \\[0.2cm]
v=c^{{l+1}\over {l+2}}_1{\varrho^{{l+1}\over {l+2}}}S\big( {{{T_*}}\phi\over {2\pi}}\big)
\end{array}
\end{equation*}
with $c_1={{2\pi(l+2)}\over {T_*}}>0$, $\varrho>0$ and $\phi\in \mathbb{T}$. We can check that  $\Psi$ is a symplectic diffeomorphism from $\mathbb{R}^+\times \mathbb{T}$ onto $\mathbb{R}^2\setminus \{0\}$. Under this transformation, Hamiltonian function (\ref{h4}) becomes
\begin{equation}\label{h7}
\begin{aligned}
H(\varrho,\phi,t)& = h(\Psi(\varrho,\phi),t)\\[0.2cm]
&= {1\over {2l+2}}c^{2l+2\over {l+2}}_1 \varrho^{2l+2\over {l+2}}-\varepsilon\sum _{j=0}^{2l} {\varepsilon^{2l-j}\over{j+1}}c^{j+1\over {l+2}}_1 {\varrho^{j+1\over {l+2}}}C^{j+1}\Big( {{{T_*}}\phi\over {2\pi}}\Big) p_j(t).
\end{aligned}
\end{equation}

After introducing two conjugate variables $\theta\in\mathbb{T}^{\mathbb{Z}}$ and $J\in \mathbb{R}^{\mathbb{Z}}$, the Hamiltonian (\ref{h7}) can
be written in the form of an autonomous Hamiltonian as follows
\begin{equation*}
H=\langle \widehat{\omega},J\rangle+{1\over {2l+2}}c^{2l+2\over {l+2}}_1 \varrho^{2l+2\over {l+2}}-\varepsilon\sum _{j=0}^{2l} {\varepsilon^{2l-j}\over{j+1}}c^{j+1\over {l+2}}_1 {\varrho^{j+1\over {l+2}}}C^{j+1}\Big( {{{T_*}}\phi\over {2\pi}}\Big) P_j(\theta),
\end{equation*}
where $P_j(\theta)$ is the shell function of the almost periodic function $p_j(t)$.

Let $[\mathfrak{C}_1,\mathfrak{C}_2]\subseteq\mathbb{R}^+$  be any bounded interval without 0, not depending on $\varepsilon$.  For any
$\varrho_0\in[\mathfrak{C}_1,\mathfrak{C}_2]$, we denote $\varrho=\varrho_0+I$ and do Taylor expansion at $\varrho_0$ for $|I|<{\varrho_0/ 2}$. Then we have
\begin{equation}\label{h8}
\begin{aligned}
H(\theta,\phi,J,I)&= \langle \widehat{\omega},J\rangle+{1\over {2l+2}}c^{2l+2\over {l+2}}_1 \varrho^{2l+2\over {l+2}}_0+{1\over {l+2}}c^{2l+2\over {l+2}}_1 \varrho^{l\over {l+2}}_0 I+O(I^2)\\[0.2cm]
&\ \ \ \ \ \ \  -\varepsilon\sum _{j=0}^{2l} {\varepsilon^{2l-j}\over{j+1}}c^{j+1\over {l+2}}_1 {(\varrho_0+I)^{j+1\over {l+2}}}C^{j+1}\Big( {{{T_*}}\phi\over {2\pi}}\Big) P_j(\theta).
\end{aligned}
\end{equation}
Denote ${\widetilde{\omega}}={1\over {l+2}}c^{2l+2\over {l+2}}_1 \varrho^{l\over {l+2}}_0$, for any $\varrho_0\in[\mathfrak{C}_1,\mathfrak{C}_2]$, we get
$${{\partial {\widetilde{\omega}}}\over {}\partial\varrho_0}={l\over {(l+2)^2}}c^{2l+2\over {l+2}}_1 \varrho^{-2\over {l+2}}_0\neq 0.$$
We therefore singer out the subsets $\mathcal{O}_{\alpha}\subseteq \mathcal{O}:=[{1\over {l+2}}c^{2l+2\over {l+2}}_1\mathfrak{C}_1^{l\over {l+2}},{1\over {l+2}}c^{2l+2\over {l+2}}_1\mathfrak{C}_2^{l\over {l+2}}]$ is the set of ${\widetilde{\omega}}$ satisfying
$${{|\langle {k},\varepsilon^l{\omega} \rangle + {\widetilde{k}}\,{\widetilde{\omega}}|} \geq {\alpha \over {\Delta\big([[(k,{\widetilde{k}})]]\big)\Delta\big(|k|+|\widetilde{k}|\big)}}},\ \ \ \ \text{for all}\ (k,{\widetilde{k}}) \in \mathbb{Z}_{\widetilde{\mathcal{S}}}^{{\mathbb{Z}}}\backslash \{\widetilde{k}\neq 0\})$$
with fixed ${{\omega}}$. From the measure estimate in Section \ref{sec:mea} with $n=1$, it follows that the set $\mathcal{O}_{\alpha}$ is a set
of positive Lebesgue measure provided that $\alpha$ is small. Let  $\mathcal{O}_h : |{\widetilde{\omega}}-\mathcal{O}_{\alpha}|<h$ denote complex neighborhoods of $\mathcal{O}_{\alpha}$.

Similar to Section \ref{sec:hal}, we can rewrite the Hamiltonian function (\ref{h8}) as
\begin{equation}\label{h9}
\begin{aligned}
H(\theta,\phi,J,I;{\widetilde{\omega}})&=N+P\\
&=e({\widetilde{\omega}})+\langle \widehat{\omega},J\rangle +{\widetilde{\omega}}\,I+P(\theta,\phi,I;{\widetilde{\omega}}),
\end{aligned}
\end{equation}
where
\begin{equation*}
\begin{aligned}
P(\theta,\phi,I;{\widetilde{\omega}})&= {1\over {2l+2}}c^{2l+2\over {l+2}}_1 \Big((\varrho_0+I)^{2l+2\over {l+2}}-\varrho^{2l+2\over {l+2}}_0-{2l+2\over {l+2}}\varrho^{l\over {l+2}}_0 I\Big)\\[0.2cm]
&\ \ \ \ \ \ \  -\varepsilon\sum _{j=0}^{2l} {\varepsilon^{2l-j}\over{j+1}}c^{j+1\over {l+2}}_1 {(\varrho_0+I)^{j+1\over {l+2}}}C^{j+1}\Big( {{{T_*}}\phi\over {2\pi}}\Big) P_j(\theta)
\end{aligned}
\end{equation*}
with $\varrho_0=g_{{\widetilde{\omega}}}({\widetilde{\omega}})$.

Hence $P$  is periodic in $\theta,\phi$  with the period $2\pi$, and real-analytic in $(\theta,\phi,I)\in \mathbb{T}^\mathbb{Z}\times\mathbb{T}\times\mathbb{R}$. Then there exist $r>0$, such that $P$ admits analytic extension in the complex neighborhood
$\{(\theta,\phi) : |\Im\ \theta|_{\infty}< r, |\Im\ \phi|< r\}$ of $\mathbb{T}^\mathbb{Z}\times\mathbb{T}$. Taking $s=\varepsilon^{1/ 2}$, then there exists $C_*>0$, depending on $l,T_*,r,h$, but not on $\varepsilon$, such that for any $|\Im\ \theta|_{\infty}< r, |\Im\ \phi|< r,  |I|< s, \widetilde{\omega}\in \mathcal{O}_h,$ we have $|||P|||_{m,r,s,h}<C_* \varepsilon$. Without losing the generality, we can assume $|||P|||_{m,r,s,h}\leq\varepsilon$, which means $s^{-1}|||P|||_{m,r,s,h}\leq \varepsilon^{1/2}.$

\subsection{The main results}
\begin{theorem}\label{thm8.1}
Every solution of (\ref{h1})  with a real analytic almost periodic function $f(t)\in AP_r({\omega})$, ${\omega}$ satisfying the nonresonance condition (\ref{b6}) is bounded. Moreover (\ref{h1}) has infinitely many almost periodic solutions.
\end{theorem}
\noindent\textbf{Proof}: If the conditions of Theorem \ref{thm8.1} hold and $0<\varepsilon<\Big({{\alpha \varepsilon_*}\over \Psi_0(\mu)\Psi_1(\rho)}\Big)^2$,\  then
$$s^{-1}|||P|||_{m,r,s,h}\leq\varepsilon^{1/2}<{{\alpha \varepsilon_*}\over \Psi_0(\mu)\Psi_1(\rho)}$$
for some $0 < \mu \leq m-w$ and $0<\rho< {r/ 2}$, where $\varepsilon_*=2^{-22}$ is an absolute positive
constant, $\Psi_0(\mu)\Psi_1(\rho)$ is defined by (\ref{f10001}), therefore the assumptions of the Theorem \ref{thm3.2} are met. Hence the existence of the invariant tori of  the Hamiltonian system (\ref{h9}) is guaranteed by Theorem \ref{thm3.2},\ the Hamiltonian system (\ref{h9}) has a real analytic invariant torus of maximal dimension and with a vector field conjugate to $(\varepsilon^l{\omega},{\widetilde{\omega}})$ for each frequency vector ${\widetilde{\omega}} \in \mathcal{O}_{\alpha}$. These families of invariant tori for all frequency vector ${\widetilde{\omega}} \in \mathcal{O}_{\alpha}$ can be visualized as invariant cylinders in the space $(t,x,\dot{x})$. These cylinders are $2\pi$-periodic in time and they become the so-called invariant tori after the identification $t\equiv t+2\pi$. Each of these tori produces a family of almost periodic solution with the frequency $(\varepsilon^l{\omega},{\widetilde{\omega}})$, all solutions with initial datum lie in the interior of some invariant cylinders, which implies that all solutions are bounded for all time.  Then system (\ref{h1}) has infinitely many almost periodic solutions  as well as the boundedness of solutions.\qed

\begin{remark}\label{rem7.2}
It follows  from the proof of Theorem \ref{thm8.1} that if the conditions of Theorem \ref{thm8.1} hold,\ then system (\ref{h1}) has infinitely many almost periodic solutions with the frequency $\{\varepsilon^l{\omega},{\widetilde{\omega}}\}$ for each frequency vector ${\widetilde{\omega}} \in \mathcal{O}_{\alpha}.$
\end{remark}

\bigskip

\section*{References}
\bibliographystyle{elsarticle-num}

\end{document}